\documentclass[12pt]{amsart}

\usepackage{tgpagella}
\usepackage{euler}
\usepackage[T1]{fontenc}
\usepackage{amsmath, amssymb}
\usepackage[hidelinks]{hyperref}
\usepackage[english]{babel}
\usepackage{mathrsfs}
\usepackage{eucal}
\usepackage[all]{xy}
\usepackage{tikz-cd}
\usetikzlibrary{decorations.pathmorphing}

\usepackage{todonotes}

\numberwithin{equation}{section}

\newtheorem{thm}{Theorem}[section]
\newtheorem*{thm*}{Theorem}
\newtheorem{lem}[thm]{Lemma}
\newtheorem{fact}[thm]{Fact}

\newtheorem{prop}[thm]{Proposition}
\newtheorem*{prop*}{Proposition}

\newtheorem{cor}[thm]{Corollary}
\newtheorem*{cor*}{Corollary}

\theoremstyle{definition}
\newtheorem{defn}[thm]{Definition}
\newtheorem*{defn*}{Definition}

\newtheorem{remark}[thm]{Remark}

\newtheorem*{question*}{Question}
\newtheorem*{Pquestion*}{Popa's question}

\newtheorem*{conv*}{Convention}

\newcommand{\C}{\mathbb{C}}
\newcommand{\N}{\mathbb{N}}
\newcommand{\Z}{\mathbb{Z}}

\newcommand{\dminus}{ 
\buildrel\textstyle\ .\over{\hbox{ 
\vrule height3pt depth0pt width0pt}{\smash-} 
}}

\def\bb{\mathbb}

\def\cal{\mathcal}

\def\M{\mathcal{M}}

\newcommand{\Flip}{\mathrm{Flip}}
\newcommand{\bseq}[1]{\mathbf{#1}}

\makeatletter

\def\dotminussym#1#2{%
  \setbox0=\hbox{$\m@th#1-$}%
  \kern.5\wd0%
  \hbox to 0pt{\hss\hbox{$\m@th#1-$}\hss}%
  \raise.6\ht0\hbox to 0pt{\hss$\m@th#1.$\hss}%
  \kern.5\wd0}
\newcommand{\dotminus}{\mathbin{\mathpalette\dotminussym{}}}

\DeclareMathOperator{\Ad}{Ad}

\DeclareMathOperator{\tr}{tr}

\newcommand\cU{{\cal U}}

\def \M{\mathcal{M}}
\def \A{\mathcal{A}}
\def \B{\mathcal{B}}
\def \two{\mathbf{2}}
\def \reg{\operatorname{reg}}

\def \boldd{\mathbf{d}}

\def \fo{\operatorname{fo}}

\allowdisplaybreaks
\makeatletter
\def\l@subsection{\@tocline{2}{0pt}{2.5pc}{5pc}{}}
\def\l@subsubsection{\@tocline{3}{0pt}{5pc}{7.5pc}{}}
\makeatother

\begin{document}

\title{Matrix ultraproducts have trivial fundamental group}

\author{Hiroshi Ando}
\address{Department of Mathematics and Informatics, Chiba University, 1-33 Yayoi-cho, Inage, Chiba
263-8522, Japan}
\email{hiroando@math.s.chiba-u.ac.jp}
\thanks{Ando was partially supported by Japan Society for the Promotion of Science (JSPS) KAKENHI 25K07024}
\author{Isaac Goldbring}
\address{Department of Mathematics\\University of California, Irvine, 340 Rowland Hall (Bldg.\# 400),
Irvine, CA 92697-3875}
\email{isaac@math.uci.edu}
\urladdr{http://www.math.uci.edu/~isaac}

\begin{abstract}
We prove that every tracial ultraproduct $\prod_\cU M_n(\bb C)$ of matrix algebras has trivial fundamental group.  In particular, this establishes that there are uncountably many nonisomorphic matrix ultraproducts, answering a question of Sorin Popa.  We actually prove the stronger result that every matrix ultraproduct has trivial first-order fundamental group and thus our techniques are model-theoretic in nature.  
\end{abstract}

\maketitle

\section{Introduction}

By a \emph{matrix ultraproduct} we mean a II$_1$ factor of the form $\prod_\cU M_n(\bb C)$, where $\cU$ is a nonprincipal ultrafilter on $\bb N$ and where the ultraproduct is the tracial ultraproduct.  Our main result concerns the fundamental group of such tracial ultraproducts, a notion first introduced by Murray and von Neumann in their seminal paper \cite{murray1943rings}, so we recall this notion first.

Suppose that $M$ is a II$_1$ factor.  For $n\in \N^{>0}$, we let $M_n(M)$ denote the algebra of $n\times n$ matrices over $M$, which is again a II$_1$ factor.  For $t\in (0,1)$, we let $M_t$ denote the corner $pMp$, where $p$ is any projection in $M$ of trace $t$; note that the isomorphism type of such a corner is independent of the choice of projection.  For arbitrary $t\in \mathbb{R}^{>0}$, choose $n\in \N$ with $n>t$ and a projection $p\in M_n(M)$ with trace $t/n$ and define $M_t:=pM_n(M)p$.  The \emph{fundamental group $\mathcal{F}(M)$ of $M$} is defined to be $\mathcal{F}(M):=\{t\in \bb R^{>0} \ : \ M_t\cong M\}$, which is readily checked to be a subgroup of the multiplicative group $\mathbb{R}^{>0}$.

Many familiar II$_1$ factors have \emph{full fundamental group}, by which we mean fundamental group equal to $\mathbb{R}^{>0}$.  For example, Murray and von Neumann proved that the hyperfinite II$_1$ factor $\mathcal{R}$ has this property; consequently, any so-called \emph{McDuff} factor, that is, a II$_1$ factor which tensorially absorbs $\mathcal{R}$, also has full fundamental group.  R\u{a}dulescu proved that
the free group factor $L(\mathbb F_\infty)$ corresponding to the free group on infinitely many generators has full fundamental group
\cite{ruadulescu1992fundamental}. Very recently, the solution by OpenAI \cite{openAI2026} to the \emph{free group factor problem} (see also Shlyakhtenko \cite{shlyakhtenko2026isomorphism}), together with work of Dykema \cite{dykema1994interpolated} and R\u{a}dulescu \cite{Radulescu94}, established that all free group factors $L(\bb F_n)$ corresponding to finite rank free groups also have full fundamental group.

Finding examples of II$_1$ factors with proper (that is, not full) fundamental group is much more challenging.  The first example known to have a proper fundamental group was obtained by Connes in \cite{connes1980factor}, where he proved that if $\Gamma$ is an ICC group with property (T), then the group factor $L(\Gamma)$ has countable fundamental group.  Building upon work of Gaboriau \cite{gaboriau2002invariants}, in \cite{popa2006class} Popa produced explicit factors with trivial fundamental group. This  deformation/rigidity approach was further developed by Popa in \cite{popa2006strong,popa2006strong2}, where it was established that any subgroup of $\mathbb{R}^{>0}$ appears as the fundamental group of a (possibly nonseparable) II$_1$ factor. 

Since matrix ultraproducts are not McDuff (in fact, they do not even have property Gamma, as shown in \cite[Theorem 5.1]{farah2014model}), it is not obvious as to whether or not they have full fundamental group.  In fact, it was an open problem to determine what exactly is the fundamental group of such  matrix ultraproducts, and whether or not the fundamental group depends on the ultrafilter.  In this paper, we settle the question by showing that all matrix ultraproducts have trivial fundamental group:

\begin{thm*}
For any nonprincipal ultrafilter $\cU$ on $\N$, we have $\mathcal{F}(\prod_\cU M_n(\bb C))=\{1\}$.
\end{thm*}

Observe now that any corner of a matrix ultraproduct is itself (isomorphic to) a matrix ultraproduct.  Indeed, given $t\in (0,1)$, let $f:\N\to \N$ be given by $f(n)=\max(1,\lfloor tn\rfloor)$.  For each $n$, pick a projection $p_n\in M_n(\bb C)$ whose rank is $f(n)$.  Set $p:=(p_n)_\cU\in \prod_\cU M_n(\C)$, a projection of trace $t$.    Since $f$ has uniformly bounded fibers, we have that $M_t\cong pMp\cong \prod_\cU M_{f(n)}(\bb C)\cong \prod_{f_*\cU}M_n(\bb C)$, where $f_*\cU$ is the \emph{push-forward} ultrafilter on $\N$ given by declaring $A\in f_*\cU$ if and only if $f^{-1}(A)\in \cU$.   

Triviality of the fundamental groups of matrix ultraproducts implies that no two corners of such ultraproducts corresponding to projections of different traces can be isomorphic.  We thus arrive at the following:  

\begin{cor*}
There are continuum many nonisomorphic matrix ultraproducts.
\end{cor*}

Let us put the previous corollary into context.  
In 2008, Sorin Popa asked Ilijas Farah whether there are uncountably
many pairwise nonisomorphic tracial matrix ultraproducts
\cite[Introduction]{farah2013model}. Farah, Hart, and Sherman showed
that there are at least two nonisomorphic matrix ultraproducts if
the continuum hypothesis fails \cite[Proposition~3.3]{farah2013model}.
Under the same hypothesis, Farah and Shelah obtained
$2^{\mathfrak c}$ pairwise nonisomorphic matrix ultraproducts, where
$\mathfrak c=2^{\aleph_0}$ \cite[Proposition~8.3]{farah2010dichotomy}.
The preceding corollary answers Popa's question in ZFC.

Since matrix ultraproducts are countably saturated and have density character the continuum, under the continuum hypothesis, isomorphism for matrix ultraproducts is equivalent to \emph{elementary equivalence}, that is, having the same first order theory.  A corollary of what we actually prove in this paper is the stronger result that there are continuum many non-elementarily equivalent matrix ultraproducts.  In particular, the first-order theories $\operatorname{Th}(M_n(\bb C))$ of matrix algebras do not converge as $n\to \infty$.  A byproduct of our proof is that there is a single $\exists\forall\exists$-sentence $\sigma$ such that $\lim_{n\to \infty}\sigma^{M_n(\bb C)}$ does not exist, in particular answering \cite[Question 6.6]{farah2014modelII} in the negative.

In order to deduce the existence of continuum many non-elementarily equivalent matrix ultraproducts, we prove that matrix ultraproducts have trivial \emph{first-order fundamental group}.  In \cite{goldbring2017theories}, the second author and Bradd Hart introduced the notion of the first-order fundamental group $\mathcal{F}_{\fo}(M)$ of a II$_1$ factor $M$, which is defined by $\mathcal{F}_{\fo}(M):=\{t\in \bb R^{>0} \ : \ M_t\equiv M\}$, where $\equiv$ denotes elementary equivalence.  In \cite{goldbring2017theories}, it was shown that $\mathcal{F}_{\fo}(M)$ is a closed subgroup of $\bb R^{>0}$ containing $\mathcal{F}(M)$ (which itself may not be closed). 

The main result of this paper is the following:

\begin{thm*}
For any nonprincipal ultrafilter $\cU$ on $\N$, we have $\mathcal{F}_{\fo}(\prod_\cU M_n(\bb C))=\{1\}$. 
\end{thm*}

The previous theorem has a number of notable consequences.  First, this theorem provides the first examples of II$_1$ factors with proper first-order fundamental group, settling \cite[Problem 48]{elayavalli202550}.  Indeed, if $N$ is a separable
elementary submodel of a matrix ultraproduct $M$, then
$N^{\mathcal V}\equiv N\equiv M$ for every nonprincipal ultrafilter
$\mathcal V$ on $\mathbb N$. Thus both $N$ and $N^{\mathcal V}$ have
trivial first-order fundamental group, and in particular
$\mathcal F(N^{\mathcal V})=\{1\}$. Second, by \cite[Proposition 4.16]{goldbring2017theories}, we see that there are continuum many distinct theories of matrix ultraproducts.  Since matrix ultraproducts are isomorphic to their opposites, two matrix ultraproducts are elementarily equivalent if and only if their unitary groups are elementarily equivalent as metric groups (see, for example, \cite{goldbring2015games}); consequently there are exactly continuum many different theories of so-called \emph{universal hyperlinear groups}, that is, metric ultraproducts of finite-dimensional unitary groups.  It follows that, under $\operatorname{CH}$, there are exactly continuum many nonisomorphic universal hyperlinear groups, complementing the result of Lupini \cite{lupini2014logic} that there are exactly $2^{\mathfrak{c}}$ many isomorphism classes of universal hyperlinear groups under the negation of $\operatorname{CH}$.  (The fact that, under $\operatorname{CH}$, there are continuum many nonisomorphic universal sofic groups was recently established by Alekseev and Thom \cite{alekseev2024non}.)  Finally, as mentioned above, the free group factors all have full fundamental group; consequently, our result offers an alternative proof that they are not elementarily equivalent to matrix ultraproducts, a result first shown by Jesse Peterson earlier this year \cite[Theorem 1.2]{peterson2026ultraproduct}.

In the next section, we give an outline of the proof of the previous theorem, highlighting the role of definable sets and groups in the proof.  The proof outline involves three major blackboxes; the proofs of these blackboxes are given in the subsequent three sections. 

To the best of our knowledge, this is the first serious application of definable groups and pseudofinite groups to the model theory of operator algebras.  In classical model theory, definable groups are central to applications as they connect model-theoretic notions such as definability, types, stability, rank, and independence with the structure of algebraic and topological groups, thereby allowing geometric model-theoretic methods to address problems in algebra, number theory, and dynamics.  A landmark example is Hrushovski’s use of the structure theory of definable groups and fields in his proof of the function-field Mordell–Lang conjecture \cite{hrushovski1996mordell}. Pseudofinite groups, typically in the form of nonprincipal ultraproducts of finite groups, are especially useful because they permit finite combinatorial problems to be studied using compactness, saturation, invariant measures, and structural results for definable sets, after which conclusions can be transferred back to sufficiently large finite groups. This perspective has produced powerful regularity and structure theorems in additive combinatorics, including results for stable subsets and sets of finite VC-dimension in finite groups \cite{conant2020group, conant2021structure}. For broader background, see the books of Poizat \cite{poizat2001stable} and Pillay \cite{pillay1996geometric}.  For some work on definable groups in continuous logic, see the papers \cite{yaacov2010stability,carmona2022definably,ivanov2021metric}.

We assume that the reader is familiar with basic tracial von Neumann algebra theory and basic continuous model theory; the relevant background in both can be found in \cite{goldbring2023model}.  The necessary facts about the first-order fundamental group and its relations to ultraproducts can be found in \cite{goldbring2017theories} and \cite{goldbring2023first}.  The proofs also involve some basic facts from the theory of representations of finite groups; these facts can be found, for example, in the textbook \cite{james2001representations}.  

After circulating a first version of our paper, we were informed by Alekseev and Thom that they had obtained the triviality of fundamental groups of matrix ultraproducts in early November, although their argument differs fairly substantially from ours; for more information on their techniques, see their recent preprint \cite{alekseev2026converge}.

We thank Jananan Arulseelan for a number of very useful comments about a first draft of the paper.

\section*{Declaration of AI Usage}

The main ideas in this paper are due to ChatGPT 6 Astra.  In late September of this year, a discussion with ChatGPT led to the proof of the existence of two non-elementarily equivalent matrix ultraproducts; subsequent discussions in October led to the full result mentioned above.  However, the proof given by ChatGPT was extremely difficult to read.  The presentation here is essentially entirely due to the authors.  The logical flow of the arguments was completely rewritten in order to make the main ideas more transparent.  More importantly, the emphasis on definability is entirely due to the authors.  In fact, ChatGPT seemed to stumble upon the notion of definability in continuous logic without actually uttering the phrase ``definable.''  We place definable groups and definable subalgebras at the forefront of the arguments to highlight their importance in the overall structure of the paper.  This is but one example of a number of existing ideas and results that were used by ChatGPT and never mentioned and/or cited, another prominent example being the Kac-Takesaki operator  or multiplicative unitary (which ChatGPT referred to by the nonstandard name ``copying unitary'') which is used to characterize group elements amongst the unitaries in the group ring.

Claude Opus 5.5 was used to produce the diagrams appearing in Section 2.

\section{Outline of the proof}

\emph{Throughout this paper, all ultrafilters are assumed to be nonprincipal}.  As stated in the introduction, the main theorem of the paper is the following:

\begin{thm}\label{thm:trivial}
For any ultrafilter $\cU$ on $\N$, we have $\mathcal{F}_{\fo}(\prod_\cU M_n(\C))=\{1\}$.
\end{thm}

We begin by making the following reduction, allowing us to focus on particular kinds of matrix ultraproducts that will be more convenient to work with in the sequel:

\begin{prop}
For any ultrafilter $\cU$, there is an ultrafilter $\mathcal{V}$ on $\N$ such that $$\mathcal{F}_{\fo}(\prod_\cU M_{n}(\C))= \mathcal{F}_{\fo}(\prod_{\mathcal{V}}M_{4^j}(\mathbb C)).$$
\end{prop}

\begin{proof}
Define $h:\N\to \N$ by $h(n)=\lfloor \log_4 n\rfloor$, whence $1\leq n/4^{h(n)}\leq 4$.  Set $s:=\lim_\cU \frac{n}{4^{h(n)}}$.  It follows that $\prod_\cU M_{n}(\C)\cong (\prod_\cU M_{4^{h(n)}}(\C))_s$.  Set $\mathcal{V}:=h_*\cU$. By {\L}o\'s's theorem, $\prod_\cU M_{4^{h(n)}}(\C)\equiv \prod_{\mathcal{V}}M_{4^j}(\C)$. As a consequence of the fact that amplification preserves elementary equivalence (see \cite[Lemma 2.5]{goldbring2023first}), the first-order fundamental group of a II$_1$ factor is equal to that of any amplification and is preserved by elementary equivalence; as a result, the proof is complete.  
\end{proof}

By the previous proposition, Theorem \ref{thm:trivial} will follow from the following theorem:

\begin{thm}\label{thm:trivialfour}
For any ultrafilter $\cU$ on $\N$, we have that $\mathcal{F}_{\fo}(\prod_\cU M_{4^j}(\C))=\{1\}$.
\end{thm}

Before moving forward, let us fix some notation.  Throughout, we let $\cU$ denote a (nonprincipal) ultrafilter on $\bb N$.  By a \emph{sequence of integers} we mean a sequence $\mathbf{d}=(d_j)_{j\in \N}$ of integers for which $\lim_\cU d_j=\infty$.  Given two sequences $\mathbf{d}$ and $\mathbf{r}$ of integers, we write $\mathbf{d}/\mathbf{r}=1$ to mean $\lim_{\cU}d_j/r_j=1$.  

For a sequence $\mathbf{d}$ of integers, we set $\M_\boldd:=\prod_\cU M_{d_j}(\C)$.  We also set $\boldd^2:=(d_j^2)_{j\in \N}$ and view $\M_{\boldd^2}$ as the ultraproduct $\prod_\cU (M_{d_j}(\C)\otimes M_{d_j}(\C))$.  We also set $\two:=(2^j)_{j\in \N}$ and $\mathbf{4}:=\mathbf{2}^2=(4^j)_{j\in \N}$, whence $\M_\mathbf{4}$ is viewed as $\prod_\cU (M_{2^j}(\C)\otimes M_{2^j}(\C))$.  

Our goal is to show that $\mathcal{F}_{\fo}(\M_{\mathbf{4}})=\{1\}$.  The majority of the work in proving that result is in proving the following:

\begin{thm}\label{thm:maintheorem}
$\mathcal{F}_{\fo}(\M_{\mathbf{4}})\subseteq 4^\Z$.
\end{thm}

We stress that the previous theorem holds for \emph{any} nonprincipal ultrafilter on $\bb N$ (even though the ultrafilter is suppressed in the notation) and this fact will be crucial to the proof of Theorem \ref{thm:trivialfour}.

In the next subsection, we spend some time outlining the proof of Theorem \ref{thm:maintheorem}.  In the subsequent subsection, we explain how to deduce Theorem \ref{thm:trivialfour} from Theorem \ref{thm:maintheorem}.

\subsection{Outline of the proof of Theorem \ref{thm:maintheorem}}\label{subsec:proofofmaintheorem}

For a sequence $\boldd$, we let $$\A_{\boldd}:=\prod_\cU (M_{d_j}(\C)\otimes 1), \qquad \B_{\boldd}:=\prod_\cU (1\otimes M_{d_j}(\C)),$$ both viewed as subalgebras of $\M_{\boldd^2}$.   For a unitary $V\in \M_{\boldd^2}$, we set $$K_{\boldd,V}:=\{u\in \mathrm{U}(\A_{\boldd}) \ : \ \|VuV^*-u\otimes u\|_2=0\}$$ and $$C_{\boldd,V}:=K_{\boldd,V}'\cap \A_\boldd.$$  Note that $K_{\boldd,V}$ is a subgroup of $\mathrm{U}(\A_\boldd)$.  In Section \ref{sec:special} below, if $(H_j)_{j\in \N}$ is a sequence of finite groups with $|H_j|=d_j$ for each $j$, then for a specifically chosen unitary $V_H\in \M_{\mathbf{d}^2}$, we show that $K_{\mathbf{d},V_H}=\prod_\cU \lambda_j(H_j)\otimes 1$, where $\lambda_j$ is the left regular representation of $H_j$.

We say that $V\in \M_{\mathbf{d}^2}$ is a \emph{special unitary} if the following three conditions are met:
\begin{enumerate}
    \item For all $u\in \mathrm{U}(\A_\boldd)$, we have $d(u,K_{\boldd,V})\leq \|VuV^*-u\otimes u\|_2$.
    \item $K_{\boldd,V}$ has exponent $2$.
    \item $C_{\boldd,V}$ is abelian.
\end{enumerate}

The following theorem will be proven in Section \ref{sec:special}:

\begin{thm}\label{thm:special}
There is a special unitary in $\M_{\boldd^2}$ if and only if there is a sequence $\mathbf{l}=(l_j)$ such that $\mathbf{d}/2^{\mathbf{l}}=1$.
\end{thm}

Here, $2^{\mathbf{l}}:=(2^{l_j})_{j\in \N}$.  In the ``if'' direction of the theorem, we will see that $V_H$ as described above is a special unitary when each $H_j\cong (\Z/2\Z)^{l_j}$ for some integer $l_j$.

Theorem \ref{thm:maintheorem} will follow from Theorem \ref{thm:special} and the following:  

\begin{thm}\label{thm:transfer}
Suppose that $\mathbf{r}$ is a sequence of integers such that $\M_\mathbf{4}\equiv \M_\mathbf{r}$.  Then there is a sequence $\mathbf{d}$ of integers such that $\mathbf{r}/\mathbf{d}^2=1$ and such that $\M_{\mathbf{d}^2}$ has a special unitary.  
\end{thm}

\begin{proof}[Proof of Theorem \ref{thm:maintheorem}]
For simplicity, set $\M:=\M_\mathbf{4}$.  Suppose that $\M\equiv \M_t$; we must show that $t\in 4^\Z$.  It suffices to consider the case that $t\in (0,1)$.  In this case, there are projections $p_j\in M_{4^j}(\C)$ with $\lim_\cU\tr_j(p_j)=t$ and $\tr_j(p_j)=r_j/4^j$ so that, setting $p:=(p_j)_\cU$, we have $\M_t=p\M p=\prod_\cU p_jM_{4^j}(\C)p_j\cong \prod_\cU M_{r_j}(\C)$.  By Theorems \ref{thm:special} and \ref{thm:transfer}, there are sequences $\mathbf{d}$ and $\mathbf{l}$ such that $\mathbf{r}/\mathbf{d}^2=1$ and $\mathbf{d}/2^{\mathbf{l}}=1$.  We thus have
$$t=\mathbf{r}/\mathbf{4}=\mathbf{d}^2/\mathbf{4}=4^{\mathbf{l}}/\mathbf{4},$$ implying that $t\in 4^{\Z}$, as desired.  
\end{proof}

In the rest of this section, we explain the proof of Theorem \ref{thm:transfer}.  First, we begin with some generalities on definability.

Fix a sequence $\mathbf{n}$ of integers and a nonnegative formula $\varphi(u,\vec v)$ whose variables range over unitaries.  Fix some parameters $\vec V\in \mathrm{U}(\M_{\mathbf{n}})$. We say that $\varphi(u,\vec V)$ is \textbf{$L$-almost near} if $d(u,Z(\varphi(\cdot, \vec V)))\leq L\varphi(u,\vec V)$ for all $u\in \mathrm{U}(\M_{\mathbf{n}})$.  This property expresses the fact that the zeroset $Z(\varphi(\cdot,\vec V))$ of $\varphi(u,\vec V)$ is a definable subset of $\mathrm{U}(\M_{\mathbf{n}})$ and the modulus of definability is $L$-Lipschitz.  (For more on definability in continuous logic, see the second author's survey \cite{goldbringspectralgap}.)  It is readily verified that being $L$-almost near is expressible by the formula
$$R_L(\varphi(\cdot,\vec V)):=\sup_{u\in \mathrm{U}(\M_{\mathbf{n}})}\inf_{v\in \mathrm{U}({\M_\mathbf{n}})}\max\left(\varphi(v,\vec V),\|u-v\|_2\dminus L\varphi(u,\vec V)\right).$$

Definability will play two crucial roles in this paper.  The first role is as a key assumption in the following theorem.  Before stating it, we need to introduce one piece of notation.

For an integer $d$, we let $\operatorname{Flip}_d$ be the unitary in $M_d(\C)\otimes M_d(\C)$ that is the flip on $\C^d\otimes \C^d$, that is, $\operatorname{Flip}_d(u\otimes v)=v\otimes u$ for all $u,v\in \bb C^d$.  For a given sequence $\mathbf{d}$ of integers, we set $\operatorname{Flip}_{\mathbf{d}}:=(\operatorname{Flip}_{d_j})_\cU\in \mathrm{U}(\M_{\mathbf{d}^2})$.  In regards to the above discussion, it is important to note that $u\otimes u=u\operatorname{Ad}(\operatorname{Flip}_{\mathbf{d}})(u)$ for all unitaries $u\in \A_{\mathbf{d}}$.

The following result will be proven in Section \ref{sec:frame}:

\begin{thm}\label{thm:frame}
Suppose that $\varphi(u,\vec V)$ is $L$-almost near and $F\in \mathrm{U}(\M_{\mathbf{n}})$ satisfies $F^2=1$.  Further suppose the following hold:
\begin{enumerate}
    \item $G:=Z(\varphi(\cdot,\vec V))$ is a subgroup of $\mathrm{U}(\M_{\mathbf{n}})$ for which $\tr(g)=0$ for all $g\in G\setminus \{1\}$.
    \item $[g,FhF^*]=0$ for all $g,h\in G$.
    \item $(G\cup FGF^*)'\cap \M_{\mathbf{n}}=\mathbb{C}$.
\end{enumerate}
Then there is a sequence $\mathbf{d}$ of integers such that $\mathbf{n}/\mathbf{d}^2=1$ and a trace-preserving $*$-isomorphism $\Phi:\M_{\mathbf{n}}\to \M_{\mathbf{d}^2}$ satisfying:
\begin{enumerate}
    \item[(i)] $\Phi(G)'\cap \M_{\mathbf{d}^2}=\B_{\mathbf{d}}$.
    \item[(ii)] $(\Phi(G)'\cap \M_{\mathbf{d}^2})'\cap \M_{\mathbf{d}^2}=\A_{\mathbf{d}}$.
    \item[(iii)] $\Phi \circ \operatorname{Ad}(F)=\operatorname{Ad}(\operatorname{Flip}_{\mathbf{d}})\circ \Phi$.
\end{enumerate}
\end{thm}

\begin{remark}
It is important to point out that, in the context of the previous theorem, we have that $\Phi(G)$ is defined by the $L$-almost near formula $\varphi(u,\Phi(\vec V))$.
\end{remark}

As in the previous remark, we will be transporting definable groups and definable subalgebras defined by almost-near formulae by isomorphisms and we will want to know that the relative commutant and relative bicommutant of such definable groups are definable and are defined by almost-near formulae in a way that is independent of the parameters.  The reason this is possible is that definable sets are precisely the sets that one can quantify over in the sense that quantification over a definable set can be expressed using an official continuous logic formula.  This is explained in great detail in \cite{goldbringspectralgap}, but we list here the special case of this phenomenon used repeatedly throughout the paper.

\begin{fact}\label{fact:quant}
Suppose that $D$ is a metric space and $p:D\to \bb R^{\geq 0}$ is a function.  Set $S:=\{x\in D \ : \ p(x)=0\}$.  Suppose that $S\not=\emptyset$ and there is $K\geq 0$ such that $d(x,S)\leq Kp(x)$ for all $x\in D$.  Then for any $L$-Lipschitz function $f:D\to [0,B]$, we have
\begin{alignat*}{2}
\sup_{x\in S}f(x)&=\sup_{x\in D} (f(x)\dminus KLp(x))\\
\inf_{x\in S}f(x)&=\inf_{x\in D}\min(B,f(x)+KLp(x)).
\end{alignat*}
\end{fact}

We now explain how to use the previous fact in the case of a definable subgroup of the unitary group.  

To explain this, we need the following averaging estimate, which is a standard
consequence of the Hilbert-space ergodic theorem
of Alaoglu and Birkhoff
\cite{AlaogluBirkhoff1940};
see also the more general von Neumann algebraic
ergodic theorem of Kov\'acs and Sz\H{u}cs
\cite[Theorem~1]{KovacsSzucs66}.
We include the short argument for completeness.

\begin{lem}\label{syn:averaging}
Let $(M,\tau)$ be a tracial von Neumann algebra, $H\leq \mathrm{U}(M)$ any subgroup, and
$B=H'\cap M$. For $x\in M_1$ (the operator-norm unit ball), one has
\[
 \|x-E_B(x)\|_2\leq\sup_{h\in H}\|[x,h]\|_2.
\]
\end{lem}
\begin{proof}
Let $\mathcal C$ be the closed convex hull in $L^2(M)$ of
$\{hxh^*:h\in H\}$. The set $\mathcal C$ is nonempty, closed, convex,
bounded, and contained in $M_1$ (since $M_1$ is closed in the $2$-norm).
It thus has a unique element $y$ of smallest $L^2$ norm. Each $\Ad h$ is an
isometry preserving $\mathcal C$, so uniqueness gives $hyh^*=y$ for all
$h\in H$ and thus $y\in B$. It follows that $E_B(hxh^*)=E_B(x)$ for all $h\in H$.  Continuity gives
$E_B(y)=E_B(x)$, and hence $y=E_B(x)$. Since $\|x-hxh^*\|_2=\|[x,h]\|_2$ for all $h\in H$, it follows that $\|x-y\|_2\leq \sup_{h\in H}\|[x,h]\|_2$, as claimed.
\end{proof}

The previous lemma immediately yields the following:

\begin{prop}\label{prop:definability}
Suppose that $\M$ is an $\aleph_1$-saturated tracial von Neumann algebra (such as a matrix ultraproduct).  Suppose that $\varphi(u,\vec v)$ is a nonnegative formula with parameters ranging over unitaries and $\vec V\in \mathrm{U}(\M)$ is such that $\varphi(u,\vec V)$ is an $L$-almost near formula defining a subgroup $G$ of $\mathrm{U}(\M)$.  Set $\B:=G'\cap \M$ and $\A:=\B'\cap \M$.  Set $$\varphi_{\operatorname{c}}(x,\vec v):=\sup_{u\in \mathrm{U}(\M)}(\|[x,u]\|_2\dminus 2L\varphi(u,\vec v))$$ and $$\varphi_{\operatorname{bc}}(x,\vec v):=\sup_{y\in \M_1}(\|[x,y]\|_2\dminus 2 \varphi_{\operatorname{c}}(y,\vec v)),$$ where the first variable in both formulae ranges over the unit ball of the algebra.  Then for all $x\in \M_1$, we have:
\begin{enumerate}
    \item 
    \begin{itemize}
        \item $\varphi_{\operatorname{c}}(x,\vec V)=\sup_{g\in G}\|[x,g]\|_2$,
        \item $d(x,\B_1)\leq \varphi_{\operatorname{c}}(x,\vec V)$, and
        \item $\B_1$ is a definable subset of $\M_1$, defined by $\varphi_{\operatorname{c}}(\cdot,\vec V)$.
    \end{itemize}  
    \item 
    \begin{itemize}
        \item $\varphi_{\operatorname{bc}}(x,\vec V)=\sup_{y\in \B_1}\|[x,y]\|_2$,
        \item $d(x,\A_1)\leq \varphi_{\operatorname{bc}}(x,\vec V)$,
        \item $\A_1$ is a definable subset of $\M_1$, defined by $\varphi_{\operatorname{bc}}(\cdot,\vec V)$.
    \end{itemize}  
\end{enumerate}
\end{prop}

Our discussion thus far has shown how to achieve definability of $\A_{\mathbf{d}}$ and $\B_{\mathbf{d}}$ \emph{modulo} the existence of an appropriate definable group.  In Section \ref{sec:affine}, we show how to arrange for the existence of such a definable group using the particular quantifier-free formula $\varphi_{\operatorname{aff}}$, namely $$\varphi_{\operatorname{aff}}(u,x_1,x_2,w,z):=\max(\|[u,x_1]\|_2,\|[u,x_2]\|_2,\|wuw^*-uzuz^*\|_2).$$

Before explaining the main result about $\varphi_{\operatorname{aff}}$ and definability, we need one piece of terminology.  Suppose that $\mathbf{d}$ and $\mathbf{r}$ are sequences of integers such that $\mathbf{d}/\mathbf{r}=1$ and such that $d_j\geq r_j$ for each $j$.  Let $p_j$ denote a projection of $M_{d_j}(\bb C)$ of trace $r_j/d_j$, so that $p_jM_{d_j}(\bb C)p_j\cong M_{r_j}(\bb C)$.  We then obtain an isomorphism $$(x_j)_\cU\mapsto (p_jx_jp_j)_\cU:\M_\mathbf{d}\to \M_{\mathbf{r}}$$ that we refer to as a \emph{corner isomorphism}. 

  The following theorem will be proven in Section \ref{sec:affine}.

\begin{thm}\label{thm:affine}
Suppose $n_j=q_j(q_j-1)$, where $q_j$ is a power of a prime such that
$\lim_{\cU}q_j=\infty$.  Then there are $b_1,b_2,W,F\in \mathrm{U}(\M_{\mathbf{n}})$ with $F^2=1$ and $L\geq 1$ such that $\varphi_{\operatorname{aff}}(u,b_1,b_2,W,F)$ is $L$-almost near and defines a subgroup $G_{\mathbf{n}}$ of $\mathrm{U}(\M_{\mathbf{n}})$ satisfying the hypotheses of Theorem \ref{thm:frame}.  Moreover, setting $d_j:=q_j-1$, there is a corner isomorphism
$$\Psi:\M_{\mathbf{n}}\to \M_{\mathbf{d}^2}$$ that satisfies:
\begin{itemize}
    \item $\Psi(G_{\mathbf{n}}'\cap \M_{\mathbf{n}})=\B_{\mathbf{d}}$ and $\Psi((G_{\mathbf{n}}'\cap \M_{\mathbf{n}})'\cap \M_{\mathbf{n}})=\A_{\mathbf{d}}$.
    \item $\Psi\circ \operatorname{Ad} F=\operatorname{Ad}(\operatorname{Flip}_{\mathbf{d}})\circ \Psi$.
\end{itemize}
\end{thm}

Now take primes $q_j>2^j$ such that $\mathbf{q}/\two=1$; the existence of such primes is a straightforward consequence of the Prime Number Theorem.  Indeed, for each $j$, let $q_j$ be the smallest prime strictly
greater than $2^j$. By the Prime Number Theorem, for every
fixed $\varepsilon>0$, we have
\[
\pi((1+\varepsilon)x)-\pi(x)
\sim \frac{\varepsilon x}{\log x}
\longrightarrow\infty
\qquad (x\to\infty).
\]
Thus, for all sufficiently large $j$, there is a prime
in $(2^j,(1+\varepsilon)2^j]$, whence
\[
1<\frac{q_j}{2^j}\leq 1+\varepsilon.
\]
Since $\varepsilon>0$ is arbitrary, it follows that
$q_j/2^j\to 1$ as $j\to\infty$, and hence
$\mathbf q/\mathbf 2=1$.

Since $\mathbf{d}/\mathbf{2}=1$ and $d_j\geq 2^j$ for all $j$, we can choose a corner isomorphism $\Lambda:\M_{\mathbf{d}^2}\to \M_\mathbf{4}$ intertwining the flip maps. 
More precisely, for each $j$, choose a projection
$e_j\in M_{d_j}(\C)$ of rank $2^j$ and set
$p_j=e_j\otimes e_j$. Since
\[
\tr(p_j)=\frac{4^j}{d_j^2}\longrightarrow 1,
\]
compression by $(p_j)_\cU$ induces a corner
isomorphism
\[
\Lambda:\M_{\mathbf{d}^2}\longrightarrow
\M_{\mathbf{4}}.
\]
By choosing the identifications
$p_jM_{d_j^2}(\C)p_j\cong M_{4^j}(\C)$
as tensor products of the corresponding
identifications
$e_jM_{d_j}(\C)e_j\cong M_{2^j}(\C)$,
we may arrange that
\[
\Lambda(\A_{\mathbf d})=\A_{\mathbf 2},
\qquad
\Lambda(\B_{\mathbf d})=\B_{\mathbf 2},
\qquad
\Lambda(\operatorname{Flip}_{\mathbf d})
=\operatorname{Flip}_{\mathbf 2}.
\]

We also have the isomorphism $\Psi:\M_{\mathbf{n}}\to \M_{\mathbf{d}^2}$ from Theorem \ref{thm:affine}. Note that $\Theta:=\Lambda\circ \Psi\colon \M_{\mathbf{n}}\to \M_{\mathbf{4}}$ satisfies  $$\Theta\circ \operatorname{Ad}(F)=\operatorname{Ad}(\operatorname{Flip}_{\mathbf{2}})\circ \Theta.$$ We have the $L$-almost near formula $\varphi_{\operatorname{aff}}(u,\vec b,W,F)$ which defines the subgroup $G_\mathbf{n}$ of $\mathrm{U}(\M_\mathbf{n})$.  Moreover, the relative bicommutant $\A:=(G_\mathbf{n}'\cap \M_{\mathbf{n}})'\cap \M_{\mathbf{n}}$ of $G_\mathbf{n}$ is defined by the formula $\varphi_{\operatorname{aff},\operatorname{bc}}(u,\vec b,W,F)$.  Applying the isomorphism $\Theta$, we have that $\varphi_{\operatorname{aff}}(u,\Theta(\vec b),\Theta(W),\operatorname{Flip}_{\mathbf{2}})$ is $L$-almost near and defines a subgroup $G_{\mathbf{4}}:=\Theta(G_{\mathbf{n}})$ of $\mathrm{U}(\M_\mathbf{4})$ whose relative bicommutant is $\A_{\mathbf{2}}=\Theta(\mathcal{A})$ and is defined by $$\varphi_{\operatorname{aff},\operatorname{bc}}(u,\Theta(\vec b),\Theta(W),\operatorname{Flip}_{\mathbf{2}}).$$  
Note that $\Theta(F)$ may differ from $\operatorname{Flip}_{\mathbf{2}}$. However, we have that $\operatorname{Ad}(\Theta(F))=\operatorname{Ad}(\operatorname{Flip}_{\mathbf{2}})$ and the formula $\varphi_{\operatorname{aff}}(u,b_1,b_2,W,F)$ depends on the last parameter $F$ only through the adjoint action $\operatorname{Ad}(F)$, which can be seen from the definition of $\varphi_{\operatorname{aff}}$. The same observation applies to the associated
formulas $\varphi_{\mathrm{aff},\operatorname{c}}$ and
$\varphi_{\mathrm{aff},\operatorname{bc}}$.

Let $V$ be a special unitary for $\M_{\mathbf{4}}$.

\begin{lem}\label{lem:specialdefinable}
Consider the formula $\psi(u,\Theta(\vec b),\Theta(W),\operatorname{Flip}_{\mathbf{2}},V)$ defined by $$\max(\varphi_{\operatorname{aff},\operatorname{bc}}(u,\Theta(\vec b),\Theta(W),\operatorname{Flip}_{\mathbf{2}}),\|VuV^*-u\operatorname{Flip}_{\mathbf{2}}u\operatorname{Flip}_{\mathbf{2}}\|_2)$$ and the formula $\psi'(x):=\max(\psi_{\operatorname{c}}(x),\varphi_{\operatorname{aff},\operatorname{bc}}(x))$ (with parameters suppressed for sake of readability).  Then:
\begin{enumerate}
    \item $\psi(u)$ is $(1+4\sqrt{2})$-almost near.
    \item $\psi'(x)$ is $\sqrt{2}$-almost near. 
\end{enumerate}
\end{lem}

\begin{proof}
For (1), suppressing parameters, let us abbreviate $\psi$ by
\[\psi(u):=\max(\delta(u),\kappa(u)),\]
where 
\[\delta(u):=\varphi_{\operatorname{aff,bc}}(u,\Theta(\vec{b}),\Theta(W),\operatorname{Flip}_{\mathbf{2}})
\]
and 
\[\kappa(u):=\|VuV^*-u\operatorname{Flip}_{\mathbf{2}}u\operatorname{Flip}_{\mathbf{2}}\|_2.\]
It is straightforward to check that $K_{\mathbf{2},V}=Z(\psi)$. Moreover, for $u\in \mathrm{U}(\M_{\mathbf{4}})$, we have $d(u,(\A_{\mathbf{2}})_1)\le \delta(u)$.
To estimate $d(u,K_{\mathbf{2},V})$, take $u\in \mathrm{U}(\M_{\mathbf{4}})$ and set $z:=E_{\A_{\mathbf{2}}}(u)\in (\A_\mathbf{2})_1$. In particular, since $z$ is a contraction, we have 
$\|u-z\|_2=d(u,(\A_{\mathbf{2}})_1).$
Since $\A_{\mathbf{2}}$ is a finite von Neumann algebra, there exists a unitary $a\in \A_{\mathbf{2}}$ such that $z=a|z|$ holds. Since $\tau\circ E_{\A_{\mathbf{2}}}=\tau$, $\tau(z^*u)=\tau(E_{\A_{\mathbf{2}}}(z^*u))=\|z\|_2^2$, so that 
\begin{align*}
    \|u-z\|_2^2&=\|u\|_2^2+\|z\|_2^2-2\operatorname{Re}\tau(z^*u)\\
    &=1-\|z\|_2^2.
\end{align*}
Therefore,
\begin{align*}
    \|z-a\|_2^2&=\tau((1-|z|)^2)\\
    &\le \tau(1-|z|^2)=1-\|z\|_2^2.
\end{align*}
Since $u-z$ is orthogonal to $\A_{\mathbf{2}}$ in $L^2(\M_{\mathbf{4}},\tau)$ and $z-a\in \A_{\mathbf{2}}$, this implies that 
\begin{align*}
    \|u-a\|_2^2&=\|u-z\|_2^2+\|z-a\|_2^2\\
    &\le 2\|u-z\|_2^2.
\end{align*}
Therefore 
\begin{align*}
    \|u-a\|_2&\le \sqrt{2}\|u-z\|_2=\sqrt{2}d(u,(\A_{\mathbf{2}})_1).
\end{align*}
We also note that $\kappa$ is $3$-Lipschitz on 
$\mathrm{U}(\M_{\mathbf{4}})$. Indeed, for $u,v\in \mathrm{U}(\M_{\mathbf{4}})$,
writing $S:=\operatorname{Flip}_{\mathbf{2}}$, we have
\[
\begin{aligned}
|\kappa(u)-\kappa(v)|
&\leq \|V(u-v)V^*\|_2
     +\|uSuS-vSvS\|_2\\
&\leq \|u-v\|_2
     +\|(u-v)SuS+vS(u-v)S\|_2\\
&\leq 3\|u-v\|_2,
\end{aligned}
\]
where the last inequality follows from the fact that
$u,v,S$ are unitaries. Moreover, since $a\in \mathrm{U}(\A_{\mathbf{2}})$, we have, by the defining property (1) of the special unitary $V$ and the fact $a\otimes a=a\operatorname{Flip}_{\mathbf{2}}a\operatorname{Flip}_{\mathbf{2}}$, that 
\[d(a,K_{\mathbf{2},V})\le \|VaV^*-a\otimes a\|_2=\kappa(a).\]
It then follows that 
\begin{align*}
d(u,Z(\psi))=d(u,K_{\mathbf{2},V})&\le \|u-a\|_2+d(a,K_{\mathbf{2},V})\\
&\le \|u-a\|_2+\kappa(a)\\
&\le \|u-a\|_2+\kappa(u)+3\|u-a\|_2\\
&\le \kappa(u)+4\sqrt{2}d(u,(\A_{\mathbf{2}})_1)\\
&\le \kappa(u)+4\sqrt{2}\delta(u)\\
&\le (1+4\sqrt{2})\psi(u).
\end{align*}
This shows that $\psi$ is $(1+4\sqrt{2})$-almost near.

For (2), by Proposition~\ref{prop:definability}, we have $\psi_{\operatorname{c}}(x)=\sup_{k\in K_{\mathbf{2},V}}\|[x,k]\|_2$, 
and hence $Z(\psi')=(K_{\mathbf{2},V}'\cap \A_{\mathbf{2}})_1=(C_{\mathbf{2},V})_1$. 
Fix $x\in (\M_{\mathbf{4}})_1$ and set $a:=E_{\A_{\mathbf{2}}}(x)$.
By Proposition~\ref{prop:definability}, we have
\[
\|x-a\|_2=d(x,\A_1)
\leq\varphi_{\operatorname{aff},\operatorname{bc}}(x).
\]
Since $K_{\mathbf{2},V}\subseteq \mathrm{U}(\A_{\mathbf{2}})$, applying
Lemma~\ref{syn:averaging} inside $\A_{\mathbf{2}}$ gives
\[
\begin{aligned}
\|a-E_{C_{\mathbf{2},V}}(a)\|_2
&\leq \sup_{k\in K_{\mathbf{2},V}}\|[a,k]\|_2\\
&=\sup_{k\in K_{\mathbf{2},V}}\|E_{\A_{\mathbf{2}}}([x,k])\|_2\\
&\leq\sup_{k\in K_{\mathbf{2},V}}\|[x,k]\|_2\\
&=\psi_{\operatorname{c}}(x),
\end{aligned}
\]
where we used the $\A_{\mathbf{2}}$-bimodularity and
the contractivity of $E_{\A_{\mathbf{2}}}$.

Since $x-a$ is orthogonal to $a-E_{C_{\mathbf{2},V}}(a)$
in $L^2(\M_{\mathbf{4}})$ and $E_{C_{\mathbf{2},V}}(a)=E_{C_{\mathbf{2},V}}(x)$,
we obtain
\[
\begin{aligned}
d(x,(C_{\mathbf{2},V})_1)^2
&=\|x-E_{C_{\mathbf{2},V}}(x)\|_2^2\\
&=\|x-a\|_2^2+\|a-E_{C_{\mathbf{2},V}}(a)\|_2^2\\
&\leq
\varphi_{\operatorname{aff},\operatorname{bc}}(x)^2
+\psi_{\operatorname{c}}(x)^2\\
&\leq 2\psi'(x)^2.
\end{aligned}
\]
Thus
\[
d(x,Z(\psi'))\leq\sqrt{2}\,\psi'(x)
\qquad(x\in (\M_{\mathbf{4}})_1),
\]
so $\psi'$ is $\sqrt{2}$-almost near.
\end{proof}

Now consider the following formulae, all of which evaluate to $0$ in $\M_{\mathbf{4}}$:

\begin{enumerate}
    \item $R_L(\varphi_{\operatorname{aff}}(\cdot,\Theta(\vec b),\Theta(W),\operatorname{Flip}_{\mathbf{2}}))$
    \item $\|\operatorname{Flip}_{\mathbf{2}}^2-1\|_2$
    \item $\sup_{u,v\in G_\mathbf{4}}\varphi_{\operatorname{aff}}(uv,\Theta(\vec b),\Theta(W),\operatorname{Flip}_{\mathbf{2}})$
    \item $\sup_{u\in G_\mathbf{4}}\varphi_{\operatorname{aff}}(u^*,\Theta(\vec b),\Theta(W),\operatorname{Flip}_{\mathbf{2}})$
    \item $\sup_{u\in G_\mathbf{4}}|\tr(u)-\tr({u})^2|$
    \item $\sup_{u,v\in G_\mathbf{4}}\|[u,{\operatorname{Flip}_{\mathbf{2}}}v{\operatorname{Flip}_{\mathbf{2}}}^*]\|_2$
    \item $\sup_{x\in {(\M_{\mathbf{4}})_1}}(\|x-\tr(x)\cdot 1\|_2\dminus \sup_{g\in G_{\mathbf{4}}}2\max(\|[x,g]\|_2,\|[x,{\operatorname{Flip}_{\mathbf{2}}}g{\operatorname{Flip}_{\mathbf{2}}}^*]\|_2))$
    \item $R_7(\psi(\cdot,\Theta(\vec b),\Theta(W),\operatorname{Flip}_{\mathbf{2}},V))$
    \item $\sup_{u\in \mathrm{U}(\A_{\mathbf{2}})}(d(u,K_{\mathbf{2},V})\dminus \kappa(u,V,\operatorname{Flip}_{\mathbf{2}}))$
    \item $\sup_{u\in K_{\mathbf{2},V}} \|u^2-1\|_2$
    \item $R_2(\psi'(\cdot,\Theta(\vec b),\Theta(W),\operatorname{Flip}_{\mathbf{2}},V))$
    \item $\sup_{x,y\in {(C_{\mathbf{2},V})_1}}\|xy-yx\|_2$
\end{enumerate}

Here, we are frequently using the unofficial formulae obtained by  quantifying over definable sets as shorthand for the official formula equivalent to them given by Fact \ref{fact:quant}. Note that (7) follows from Lemma \ref{syn:averaging}, since $G_{\mathbf4}$ and its flip conjugate commute
and have trivial joint relative commutant.  Note also that in (9), the predicate $d(u,K_{\mathbf{2},V})$ is indeed definable, meaning that it can be equivalently expressed using a uniform limit of official formulae and this predicate is definable uniformly from the data defining $K_{\mathbf{2},V}$.   Finally, note that we are being a bit abusive with the notation in item (11) as the predicate $R_2$ was initially only defined for formulae with object variables ranging over unitaries.  The formula in (11) is meant to have variables range over the unit ball; the adaptation to this context is straightforward.

Note that, if rewritten in terms of official formulae, the formula obtained by taking the maximum of the above twelve formulae, would be $\forall_5$.  In the appendix, we show how to bring the quantifier complexity down to $\forall_2$.

We can now prove Theorem \ref{thm:transfer}.  Let $\sigma$ denote the sentence obtained by taking the maximum of the numbered formulae and existentially quantifying over the parameters.  Then $\sigma^{\M_{\mathbf{4}}}=0$, whence $\sigma^{\M_{\mathbf{r}}}=0$.  Since $\M_\mathbf{r}$ is countably saturated, we can witness the quantifiers by some parameters $\tilde{b}_1,\tilde{b}_2,\tilde{W},\tilde{F},\tilde{V}\in \mathrm{U}(\M_{\mathbf{r}})$.  Then $\varphi_{\operatorname{aff}}(u,\tilde{b}_1,\tilde{b}_2,\tilde{W},\tilde{F})$ is an $L$-almost near formula that defines a group $G_{\mathbf{r}}$ satisfying the hypotheses of Theorem \ref{thm:frame}, which implies that there is a sequence $\mathbf{t}$ with $\mathbf{r}/\mathbf{t}^2=1$ and an isomorphism $\Pi:\M_{\mathbf{r}}\to \M_{\mathbf{t}^2}$ satisfying the conclusion of Theorem \ref{thm:frame}.  (Since $\mathbf{d}$ was already taken to be $\mathbf{q}-1$, we rename the desired sequence $\mathbf{t}$.)  Using the rest of the conjuncts (and Proposition \ref{prop:definability} to show that these conjuncts do indeed define $K_{\mathbf{t},\Pi(\tilde{V})}$ and $C_{\mathbf{t},\Pi(\tilde{V})}$), it follows that $\M_{\mathbf{t}^2}$ has a special unitary $\Pi(\tilde{V})$, finishing the proof of Theorem \ref{thm:transfer}.

\begin{figure}[tbp]
\centering
\[
\begin{tikzcd}[column sep=4em, row sep=3em]
\M_{\bseq n}
  \arrow[r,"\Psi"] \arrow[d,"\Ad F"']
  \arrow[rr, bend left=30, "\Theta"]
& \M_{\bseq d^2}
  \arrow[r,"\Lambda"] \arrow[d,"\Ad\Flip_{\bseq d}"]
& \M_{\bseq 4}
  \arrow[d,"\Ad\Flip_{\bseq 2}"]
& \M_{\bseq r}
  \arrow[r,"\Pi"] \arrow[d,"\Ad\widetilde F"']
& \M_{\bseq t^2}
  \arrow[d,"\Ad\Flip_{\bseq t}"]
\\
\M_{\bseq n}
  \arrow[r,"\Psi"']
  \arrow[rr, bend right=30, "\Theta"']
& \M_{\bseq d^2}
  \arrow[r,"\Lambda"']
& \M_{\bseq 4}
& \M_{\bseq r}
  \arrow[r,"\Pi"']
& \M_{\bseq t^2}
\end{tikzcd}
\]
\caption{The left squares are Theorem \ref{thm:affine} and $\Lambda(\mathrm{Flip}_{\mathbf d})=\mathrm{Flip}_{\mathbf 2}$;
the right square is conclusion~(iii) of Theorem \ref{thm:frame}.}
\label{fig:thm26-flips}
\end{figure}
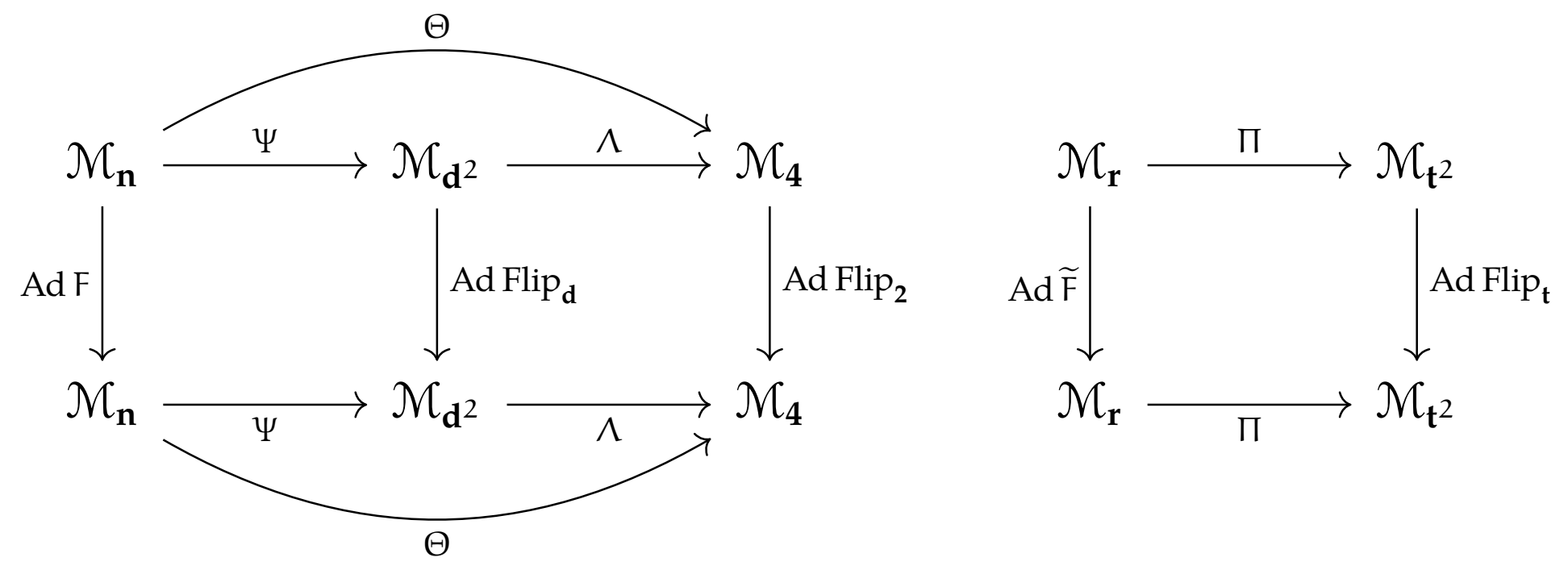

Since there are a lot of moving parts in the proof, we summarize it in Figures~\ref{fig:thm26-objects} and~\ref{fig:thm26-flips}.
\begin{figure}[p]
\centering
\[
\begin{tikzcd}[column sep=2.8em, row sep=1.9em]
G_{\bseq n}'\cap \M_{\bseq n}
  \arrow[r,"\Psi"] \arrow[d,hook]
& \B_{\bseq d}
  \arrow[r,"\Lambda"] \arrow[d,hook]
& \B_{\bseq 2}
  \arrow[d,hook]
\\
\M_{\bseq n}
  \arrow[r,"\Psi","\cong"']
& \M_{\bseq d^2}
  \arrow[r,"\Lambda","\cong"']
& \M_{\bseq 4}
  \arrow[dddddd, leftrightsquigarrow, out=0, in=0, looseness=1.1,
         "{\equiv}"]
\\
\A
  \arrow[r,"\Psi"] \arrow[u,hook']
& \A_{\bseq d}
  \arrow[r,"\Lambda"] \arrow[u,hook']
& \A_{\bseq 2}
  \arrow[u,hook']
\\
G_{\bseq n}
  \arrow[r,"\Psi"] \arrow[u,hook']
& \Psi(G_{\bseq n})
  \arrow[r,"\Lambda"] \arrow[u,hook']
& G_{\bseq 4}
  \arrow[u,hook']
\\
(\vec b,W,F)
  \arrow[rr,"\Theta"]
&
& (\Theta(\vec b),\Theta(W),\Flip_{\bseq 2},V)
\\
&
& \begin{gathered}K_{\bseq 2,V}\subseteq \mathrm{U}(\A_{\bseq 2})\\
    C_{\bseq 2,V}=K_{\bseq 2,V}'\cap \A_{\bseq 2}\end{gathered}
\\[1.5em]
& \B_{\bseq t}
  \arrow[d,hook]
& G_{\bseq r}'\cap \M_{\bseq r}
  \arrow[l,"\Pi"'] \arrow[d,hook]
\\
& \M_{\bseq t^2}
& \M_{\bseq r}
  \arrow[l,"\Pi"',"\cong"]
\\
& \A_{\bseq t}
  \arrow[u,hook']
& (G_{\bseq r}'\cap \M_{\bseq r})'\cap \M_{\bseq r}
  \arrow[l,"\Pi"'] \arrow[u,hook']
\\
& \Pi(G_{\bseq r})
  \arrow[u,hook']
& G_{\bseq r}
  \arrow[l,"\Pi"'] \arrow[u,hook']
\\
& \Pi(\widetilde V)
& (\tilde b_1,\tilde b_2,\widetilde W,\widetilde F,\widetilde V)
  \arrow[l,"\Pi"']
\\
& \begin{gathered}K_{\bseq t,\Pi(\widetilde V)}\subseteq \mathrm{U}(\A_{\bseq t})\\
    C_{\bseq t,\Pi(\widetilde V)}=K_{\bseq t,\Pi(\widetilde V)}'\cap \A_{\bseq t}\end{gathered}
&
\end{tikzcd}
\]
\caption{The ultraproducts and isomorphisms in the proof of Theorem \ref{thm:transfer}.}
\label{fig:thm26-objects}
\end{figure}
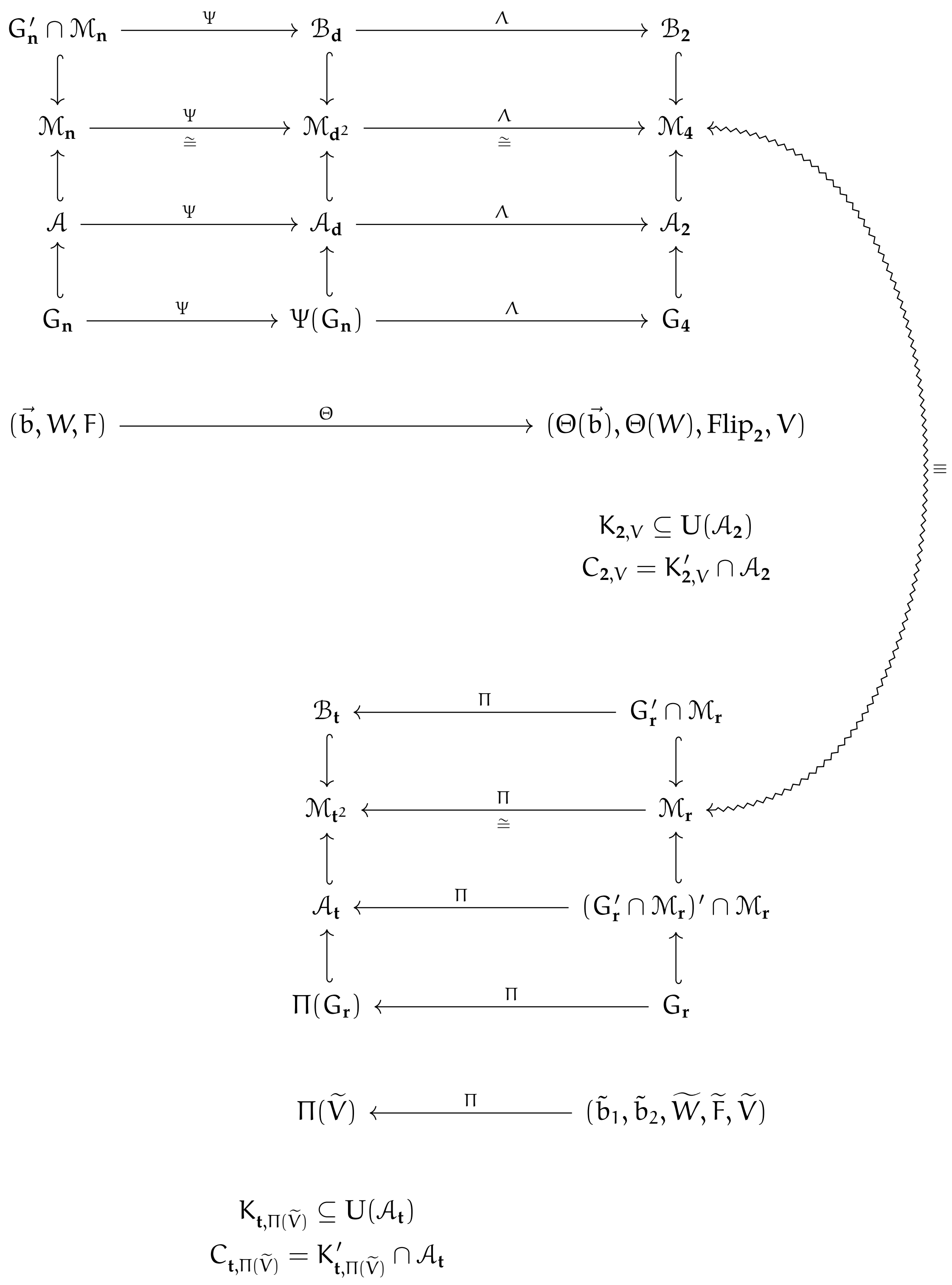

At this point, we have established Theorem \ref{thm:maintheorem} using three blackboxes, namely Theorem \ref{thm:frame}, Theorem \ref{thm:affine}, and Theorem \ref{thm:special}.  The following three sections are devoted to the proofs of these three theorems in that order.

\subsection{Finishing the proof of Theorem \ref{thm:trivialfour}}
Let us return to the task of deducing Theorem \ref{thm:trivialfour} from Theorem \ref{thm:maintheorem}.  Given a sequence $\mathbf{n}$, let $\sqrt{\mathbf{n}}$ denote the sequence $(\lfloor \sqrt{n_j}\rfloor)_{j\in \N}$.  To deduce Theorem \ref{thm:trivialfour}, it will suffice to prove:

\begin{prop}\label{prop:sqrt}
If $t\in \mathcal{F}_{\operatorname{fo}}(\M_{\mathbf{n}})$, then $\sqrt{t}\in \mathcal{F}_{\operatorname{fo}}(\M_{\sqrt{\mathbf{n}}})$.
\end{prop}

Assuming Proposition \ref{prop:sqrt}, we can now finish the proof of Theorem \ref{thm:trivialfour}.  Since we will now be dealing with many different ultrafilters, we abandon our earlier convention and make the ultrafilter explicit in the notation.  More precisely, for an ultrafilter $\cU$ and a sequence $\mathbf{n}=(n_j)_{j\in \N}$ of integers (which we still assume satisfies $\lim_\cU n_j=\infty$), we let $\M_{\mathbf{n},\cU}$ denote $\prod_\cU M_{n_j}(\C)$. 

\begin{proof}[Proof of Theorem \ref{thm:trivialfour}]
Let $f:\N\to \N$ be $f(j):=\lfloor j/2\rfloor$.  We first note the following consequence of Proposition \ref{prop:sqrt}, namely:  if $t\in \mathcal{F}_{\operatorname{fo}}(\M_{\mathbf{4},\cU})$, then $\sqrt{t}\in \mathcal{F}_{\operatorname{fo}}(\M_{\mathbf{4},f_*\cU})$.  Indeed, Proposition \ref{prop:sqrt} implies that $\sqrt{t}\in \mathcal{F}_{\operatorname{fo}}(\M_{\mathbf{2},\cU})$.  It remains to observe that $\M_{\mathbf{2},\cU}$ is either isomorphic to $\M_{\mathbf{4},f_*\cU}$ or $M_2(\M_{\mathbf{4},f_*\cU})$, depending on whether or not $\cU$ gives measure $1$ to the even numbers or the odd numbers.  However, both of these matrix ultraproducts have the same first-order fundamental group, completing the proof of the observation.

Inductively define a sequence $(\cU_i)_{i\in \N}$ of ultrafilters on $\N$ by $\cU_0:=\cU$ and $\cU_{i+1}:=f_*(\cU_i)$.  Suppose that $t\in \mathcal{F}_{\operatorname{fo}}(\M_{\mathbf{4},\cU})$.  Recursively applying the observation of the previous paragraph yields $t^{2^{-i}}\in \mathcal{F}_{\operatorname{fo}}(\M_{\mathbf{4},\cU_i})$.  Since Theorem \ref{thm:maintheorem} holds for \emph{any} ultrafilter, we have that $\mathcal{F}_{\operatorname{fo}}(\M_{\mathbf{4},\cU_i})\subseteq 4^\Z$.  It follows that $t^{2^{-i}}=4^{k_i}$ for some integer $k_i$ and thus $t=4^{2^ik_i}$.  Since this holds for every $i\geq 0$, it follows that $t=1$.    
\end{proof}

It now remains to prove Proposition \ref{prop:sqrt}.  This proposition will follow from the following more general result:

\begin{prop}\label{prop:sqrtpreserve}
Suppose that $\mathbf{n}$ and $\mathbf{m}$ are sequences and $\cU$ and $\mathcal{V}$ are ultrafilters for which $\M_{\mathbf{n},\cU}\equiv \M_{\mathbf{m},\mathcal{V}}$.  Then $\M_{\sqrt{\mathbf{n}},\cU}\equiv \M_{\sqrt{\mathbf{m}},\mathcal{V}}$.
\end{prop}

To see that Proposition \ref{prop:sqrtpreserve} implies Proposition \ref{prop:sqrt}, suppose $t\in \mathcal{F}_{\operatorname{fo}}(\M_{\mathbf{n},\cU})$ and set $m_j:=\max(1,\lfloor tn_j\rfloor)$.  Note that $\M_{\mathbf{m},\cU}\cong (\M_{\mathbf{n},\cU})_t$; the latter algebra is elementarily equivalent to $\M_{\mathbf{n},\cU}$ by assumption.  Proposition \ref{prop:sqrtpreserve} implies that $\M_{\sqrt{\mathbf{n}},\cU}\equiv \M_{\sqrt{\mathbf{m}},\mathcal{U}}$.  Since $\sqrt{\mathbf{m}}/\sqrt{\mathbf{n}}=\sqrt{t}$, we have that $\M_{\sqrt{\mathbf{m}},\cU}\cong (\M_{\sqrt{\mathbf{n}},\mathcal{U}})_{\sqrt{t}}$, which yields the desired result.

It thus remains to prove Proposition \ref{prop:sqrtpreserve}.  We need one lemma.  To simplify the notation, when discussing the formula $\varphi_{\operatorname{aff}}$ from the previous subsection, we lump all of its parameter variables into one tuple.

\begin{lem}\label{lem:affinesqrt}
For any ultrafilter $\cU$ and sequence $\mathbf{n}$, there are parameters $\vec z\in \mathrm{U}(\M_{\mathbf{n},\cU})$ such that the zeroset of $\varphi_{\operatorname{aff}}(\cdot,\vec z)$ in $\M_{\mathbf{n},\cU}$ defines a subgroup $G$ of $\mathrm{U}(\M_{\mathbf{n},\cU})$ satisfying the hypotheses of Theorem \ref{thm:frame} and for which, setting $\mathcal{A}$ for the relative bicommutant of $G$ in $\M_{\mathbf{n},\cU}$, we have $\mathcal{A}\cong \M_{\sqrt{\mathbf{n}},\cU}$.
\end{lem}

\begin{proof}
Take prime numbers $q_j>\max(3,\sqrt{n_j})$ and such that $\lim_\cU q_j/\sqrt{n_j}=1$; once again, this is possible by the Prime Number Theorem.  Set $N_j:=q_j(q_j-1)$, whence $\lim_\cU N_j/n_j=1$.  We then have unitary parameters $\vec w$ in $\mathrm{U}(\prod_\cU M_{N_j}(\C))$ for which the zeroset of $\varphi_{\operatorname{aff}}(\cdot, \vec w)$ defines a subgroup $G^\#$ of $\mathrm{U}(\M_{\mathbf{N},\cU})$ satisfying the hypotheses of Theorem \ref{thm:frame} and for which, setting $\mathcal{A}^\#$ for the relative bicommutant of $G^\#$ in $\mathcal{M}_{\mathbf{N},\cU}$, we have that $\mathcal{A}^\#\cong \M_{\mathbf{d},\cU}$, where $d_j:=(q_j-1)$.  Since $\M_{\mathbf{n},\cU}\cong \M_{\mathbf{N},\cU}$ and $\M_{\mathbf{d},\cU}\cong \M_{\sqrt{\mathbf{n}},\cU}$, the proof is complete.
\end{proof}

\begin{proof}[Proof of Proposition \ref{prop:sqrtpreserve}]
Take a unitary tuple $\vec z\in \mathrm{U}(\M_{\mathbf{n},\cU})$ as in Lemma \ref{lem:affinesqrt}.  By countable saturation of $\M_{\mathbf{m},\mathcal{V}}$, there is a tuple $\vec w\in \mathrm{U}(\M_{\mathbf{m},\mathcal{V}})$ such that $\operatorname{tp}^{\M_{\mathbf{n},\cU}}(\vec z)=\operatorname{tp}^{\M_{\mathbf{m},\mathcal{V}}}(\vec w)$.  
By the construction in Lemma~\ref{lem:affinesqrt}, the formula
$\varphi_{\operatorname{aff}}(\cdot,\vec z)$ is $L$-almost near, and
the analogues of conditions~(1)--(7) in
Subsection~\ref{subsec:proofofmaintheorem} vanish at $\vec z$.
Since $\vec z$ and $\vec w$ have the same type, these conditions also
vanish at $\vec w$. Thus $\varphi_{\operatorname{aff}}(\cdot,\vec w)$
is $L$-almost near and its zeroset $G^\#$ satisfies the hypotheses of
Theorem~\ref{thm:frame}.

Let $\mathcal{A}^\#$ be the relative bicommutant of $G^\#$ in
$\M_{\mathbf{m},\mathcal{V}}$. By Proposition~\ref{prop:definability},
the formula $\varphi_{\operatorname{aff},\operatorname{bc}}$, with
parameters $\vec z$ and $\vec w$, defines the unit balls of
$\mathcal{A}$ and $\mathcal{A}^\#$, respectively, with the same linear
distance bound. By relativization to definable sets
\cite[Theorem~2.13]{goldbringspectralgap}, every sentence about these
bicommutants can be evaluated by the same definable predicate in the
ambient algebras, at $\vec z$ and $\vec w$, respectively.
Equality of types therefore gives $\mathcal{A}\equiv\mathcal{A}^\#$.
By Theorem~\ref{thm:frame}, there is a sequence $\mathbf d$ with
$\lim_{\mathcal V}m_j/d_j^2=1$ such that
$\mathcal A^\#\cong \M_{\mathbf d,\mathcal V}
\cong \M_{\sqrt{\mathbf m},\mathcal V}$.
Since $\mathcal A\cong\M_{\sqrt{\mathbf n},\cU}$ by
Lemma~\ref{lem:affinesqrt}, it follows that
$\M_{\sqrt{\mathbf n},\cU}\equiv\M_{\sqrt{\mathbf m},\mathcal V}$.

\end{proof}

The proof of Theorem \ref{thm:trivialfour}, and thus that of Theorem \ref{thm:trivial}, has now been established (modulo the three blackboxes to be proven in the following three sections).

\subsection{A couple of preliminary lemmas}The following standard fact is recalled in
\cite[Section~1.6]{Popa2014}.
We include the short proof for completeness.

\begin{lem}[{\cite[Section~1.6]{Popa2014}}]
\label{rec:internal-algebra}
Let $B_j\subseteq M_{n_j}(\C)$ be unital $*$-subalgebras and let
$B=\prod_{\cU} B_j\subseteq \M_{\mathbf{n}}$.  Then
\[
 B'\cap \M_{\mathbf{n}}=\prod_{\cU}\bigl(B_j'\cap M_{n_j}(\C)\bigr).
\]
\end{lem}

\begin{proof}
The inclusion $\supseteq$ is clear.  For
the other inclusion, take $x=(x_j)_\cU\in B'\cap \M_{\mathbf{n}}$ and set
\[
 a_j=\sup_{v\in\mathrm U(B_j)}\|[x_j,v]\|_{2}.
\]
Note that $\lim_\cU a_j=0$.  We have a trace-preserving conditional expectation $E_j:M_{n_j}(\C)\to B_j'\cap M_{n_j}(\C)$ defined by 
\[
 E_j(x_j)=\int_{\mathrm U(B_j)}v x_jv^*\,dv
 \in B_j'\cap M_{n_j}(\C).
\]
Since
$\|x_j-E_j(x_j)\|_2\leq a_j$, we have that $x=(E_j(x_j))_\cU\in \prod_\cU (B_j'\cap M_{n_j}(\C))$,
as required.
\end{proof}
The same averaging argument gives the following group version.
\begin{lem}
\label{rec:internal-group}
Let $H_j$ be finite groups and
$\rho_j:H_j\to\mathrm U(n_j)$ be unitary representations.  Define
the subgroup $H:=\prod_\cU \rho_j(H_j)$ of the unitary group of $\M_{\mathbf{n}}$.
We then have
\[
 H'\cap \M_{\mathbf{n}}=\prod_{\cU}(\rho_j(H_j)'\cap M_{n_j}(\C))
\]
and 
\[
 (H'\cap \M_{\mathbf{n}})'\cap \M_{\mathbf{n}}=\prod_{\cU}\left((\rho_j(H_j)'\cap M_{n_j}(\C))'\cap M_{n_j}(\C)\right).
\]
\end{lem}

\begin{proof}
For the commutant statement, note that the inclusion $\supseteq$ is clear.  For the reverse inclusion, given $x_j\in M_{n_j}(\C)$, set 

\[
 E_j(x_j)=\frac1{|H_j|}\sum_{h\in H_j}
       \rho_j(h)x_j\rho_j(h)^*
       \]
and note that $E_j(x_j)$ belongs to $\rho_j(H_j)'\cap M_{n_j}(\C)$.
The inequality
    
\[
 \|x_j-E_j(x_j)\|_2
 \leq\sup_{h\in H_j}\|[x_j,\rho_j(h)]\|_2
\]
then establishes the other inclusion.  
The bicommutant identity then follows from 
Lemma~\ref{rec:internal-algebra}.
\end{proof}

\section{Proof of Theorem \ref{thm:frame}}\label{sec:frame}

In this section, we work towards the proof of Theorem \ref{thm:frame}.  The first goal is to show that every definable subgroup of the unitary group of $\M_\mathbf{n}$ (satisfying an extra technical condition) is, after passing to an isomorphic matrix ultraproduct, a pseudofinite group, that is, an ultraproduct of finite groups.  We first prove something slightly weaker than this, for which we need a definition.  From now on, we simplify notation and write $\mathrm{U}(n)$ instead of $\mathrm{U}(M_n(\bb C))$ (and similarly for other finite-dimensional unitary groups).

\begin{defn}
Suppose that $G$ is a group, $n\in \N$, and $\epsilon>0$.  Then a function $f:G\to \mathrm{U}(n)$ is called an \emph{$\epsilon$-almost representation} if $f(e)=1$ and $$\sup_{g,h\in G}\|f(gh)-f(g)f(h)\|_2\leq \epsilon.$$
\end{defn}

\begin{thm}\label{thm:pseudofinite}
Suppose that $\varphi(u,\vec V)$ is an $L$-almost-near formula with parameters in $\M_{\mathbf{n}}$ that defines a subgroup $G$ of $\mathrm{U}(\M_{\mathbf{n}})$.  Further suppose that there is $s>0$ such that $d(u,v)>s$ for all pairs of distinct elements $u,v\in G$.  Then there are finite groups $G_j$, real numbers $\epsilon_j>0$ with $\lim_\cU \epsilon_j=0$, and $\epsilon_j$-almost representations $f_j:G_j\to \mathrm{U}(n_j)$ such that $G=\prod_\cU f_j(G_j)$. 
\end{thm}

\begin{proof}
Write $\vec V=(\vec V_j)_\cU$ and set $$R_{L,j}(\varphi(\cdot,\vec V_j)):=\sup_{u\in \mathrm{U}(n_j)}\inf_{v\in \mathrm{U}(n_j)}\max(\varphi(v,\vec V_j),\|u-v\|_2\dminus L\varphi(u,\vec V_j)).$$ 

Set $\varphi_j(u,\vec V_j):=\varphi^{M_{n_j}(\C)}(u,\vec V_j)$, and 
\[
 r_j=R_{L,j}+\varphi_j(1,\vec V_j)+\frac1j,
 \qquad
 S_j=\{u\in\mathrm{U}(n_j):\varphi_j(u,\vec V_j)\leq r_j\}.
\]
Since $1\in G$, we have $\lim_\cU r_j=0$ and $1\in S_j$.  For every 
$u\in\mathrm U(n_j)$, choose $T_j(u)\in S_j$ satisfying
\begin{equation}\label{rec:rounding}
 \|u-T_j(u)\|_2\leq L\varphi_j(u)+r_j.
\end{equation}

Define
\begin{align*}
 e_j&=\sup_{u,v\in S_j}
       \min\{\|u-v\|_2,s\dminus \|u-v\|_2\},\\
 b_j&=\max\left\{
       \sup_{u,v\in S_j}\varphi_j(uv),
       \sup_{u\in S_j}\varphi_j(u^*)\right\},
       \qquad t_j=Lb_j+r_j.
\end{align*}
It is clear that $\lim_\cU e_j=\lim_\cU b_j=\lim_\cU t_j=0$.  Without loss of generality, we may assume that $e_j,t_j<s/16$ for all $j\in \N$.  For $u,v\in S_j$,
the definition of $e_j$ gives the dichotomy
\begin{equation}\label{rec:dichotomy}
 \|u-v\|_2\leq e_j
 \quad\hbox{or}\quad
 \|u-v\|_2\geq s-e_j.
\end{equation}
Define a relation $\sim_j$ on $S_j$ by declaring $u\sim_jv$ when $\|u-v\|_2<s/2$.  Using \eqref{rec:dichotomy}, it is easy to verify that $\sim_j$ is an equivalence
relation on $S_j$.  Moreover, each equivalence class has
diameter at most $e_j$ and different classes are separated by
at least $s-e_j$.  There are finitely many classes since $S_j$ is totally bounded.  Set $G_j:=S_j/\sim_j$ and write $[u]_j$ for the $\sim_j$-class of $u\in S_j$.  We turn $G_j$ into a group by defining
\[
 [u]_j[v]_j:=[T_j(uv)]_j.
\]
To see that multiplication is well-defined, for $u,v\in S_j$, \eqref{rec:rounding} gives
$\|T_j(uv)-uv\|_2\leq t_j$.  Consequently, if $u\sim_ju'$ and $v\sim_jv'$,
then
\[
 \|T_j(uv)-T_j(u'v')\|_2
 \leq2t_j+\|uv-u'v'\|_2
 \leq2t_j+2e_j<s/2.
\]
To check associativity,
the representatives
$T_j(T_j(uv)w)$ and $T_j(uT_j(vw))$ are each at distance at most
$2t_j$ from $uvw$.  Their mutual distance is at most
$4t_j<s/2$, so they have the same class.  The class of $1$ is an identity since
$\|T_j(u)-u\|_2\leq t_j$.  Finally, the class $[T_j(u^*)]_j$ is a two-sided
inverse of $[u]_j$ as the representatives of their two products
are at distance at most $2t_j<s/2$ from $1$.  We have thus established that $G_j$ is
a finite group.

To define the maps $f_j:G_j\to \mathrm{U}(n_j)$, we simply choose a representative $f_j(g)\in S_j$ of each class, making sure to choose
$f_j([1]_j)=1$.  The representative $f_j(gh)$ and
$T_j(f_j(g)f_j(h))$ are in the same class and hence are within
$e_j$.  Consequently,
\begin{equation}\label{rec:explicit-product}
 \sup_{g,h\in G_j}\|f_j(gh)-f_j(g)f_j(h)\|_2
 \leq e_j+t_j
\end{equation}
and so each $f_j$ is an $\epsilon_j$-almost representation with $\epsilon_j:=e_j+t_j$.

Since $\varphi_j(f_j(g_j),\vec V_j)\leq r_j$, it follows that $\prod_\cU f_j(G_j)\subseteq G$.  Conversely, consider
$u=(u_j)_\cU\in G$.  The point
$T_j(u_j)$ belongs to some class $g_j\in G_j$, and
\[
 \lim_\cU\|u_j-f_j(g_j)\|_2
 \leq \lim_\cU(L\varphi_j(u_j,\vec V_j)+r_j+e_j)=0,
\]
whence $u=(f_j(g_j))_\cU$, finishing the proof.
\end{proof}

\begin{remark}\label{rem:generalpredicate}
Note that in the proof of Theorem \ref{thm:pseudofinite}, we obtained the following inequality (on a $\mathcal{U}$-large set):
\begin{equation}\|f_j(g)-f_j(h)\|_2\ge s-e_j,\qquad g,h\in G_j,\,g\neq h.\label{eq:fjestimate}
\end{equation}
This inequality follows from the fact that different $\sim_j$ classes are separated by at least $s-e_j$ in the 2-norm. This inequality will be used later to obtain (\ref{eq:rho-separation}). 

Note also that we did not actually use the fact that $\varphi$ was a formula of continuous logic.  All that was actually used in the proof is that there was a sequence of predicates $\varphi_j:\mathrm{U}(n_j)\to [0,C]$ with a common modulus of uniform continuity for which the predicate $\varphi:\mathrm{U}(\M_{\mathbf{n}})\to [0,C]$ given by $\varphi((u_j)_\cU):=\lim_\cU \varphi_j(u_j)$ is $L$-almost near (in the same sense as was used for actual formulae).  (Note that the proof of Theorem \ref{thm:pseudofinite} in this more general context establishes $\lim_\cU R_{L,j}=0$ not from \L os' theorem but from the $L$-almost nearness of $\varphi$.)  This more general statement of the theorem will be used in Section \ref{sec:special}.
\end{remark}

We will need the following fact, which is known as \emph{uniform flexible stability}, and roughly states that almost representations are near actual representations (after a small amplification).  The result is due to Gowers and Hatami and we refer to \cite[Theorem 1.6]{de2019operator} for a reference.  Since the statement is slightly different from the published version, we offer a brief proof sketch.

\begin{fact}\label{fact:flexible}
There are functions $a(\epsilon),b(\epsilon)\in (0,1)$ defined for sufficiently small $\varepsilon>0$
such that $\lim_{\epsilon\to 0}a(\epsilon)=\lim_{\epsilon\to 0}b(\epsilon)=0$ satisfying the following property:  whenever $H$ is a finite group and $f:H\to \mathrm{U}(n)$ is an $\epsilon$-almost representation, then there is an integer $D$ such that $n\leq D\leq (1+a(\epsilon))n$ and a unitary representation $\rho:H\to \mathrm{U}(D)$ such that $$\sup_{g\in H}\|\rho(g)-\hat f(g)\|_{2,D}\leq b(\epsilon),$$ where $J:\C^n\to \C^D$ is an isometry and $\hat f:H\to \mathrm{U}(D)$ is given by $\hat f(g)=Jf(g)J^*+(1-JJ^*)$.
\end{fact}

\begin{proof}[Proof Sketch]
\cite[Theorem 1.6]{de2019operator} provides $D$, $\rho$ and $J$ as in the statement of the fact with $(D-n)/D=O(\epsilon^2)$ and
$\sup_{g\in H}\|f(g)-J^*\rho(g)J\|_{2,n}=O(\epsilon)$. Setting $p=JJ^*$, we have
$p(\rho(g)-\hat f(g))p=J(J^*\rho(g)J-f(g))J^*$, while the remainder
$(1-p)(\rho(g)-\hat f(g))+p(\rho(g)-\hat f(g))(1-p)$ is a sum of two
operators of rank at most $D-n$ and norm at most $2$, so
$\|\rho(g)-\hat f(g)\|_2\le\|J^*\rho(g)J-f(g)\|_2+4\sqrt{(D-n)/D}=O(\epsilon)$.
\end{proof}
Combining Theorem \ref{thm:pseudofinite} with Fact \ref{fact:flexible}, we have:

\begin{cor}\label{cor:pseudofinite}
Suppose that $\varphi(u)$ is an $L$-almost
near formula for some $L>0$, with parameters
in $\M_{\mathbf n}$, whose zeroset
$G:=Z(\varphi)$ is a subgroup of
$\mathrm U(\M_{\mathbf n})$.
Further suppose that there is $s>0$ such that $d(u,v)>s$ for all pairs of distinct elements $u,v\in G$.  Then there is a sequence $\mathbf{D}$ with $D_j\geq n_j$ for all $j$ and such that $\mathbf{D}/\mathbf{n}=1$, a sequence of finite groups $G_j$, and a sequence of unitary representations $\rho_j:G_j\to \mathrm{U}(D_j)$ such that, letting $\Phi:\M_{\mathbf{n}}\to \M_{\mathbf{D}}$ be the inverse of the natural corner isomorphism, we have $\Phi(G)=\prod_\cU \rho_j(G_j)$.
Moreover, for $\mathcal U$-almost all $j$,
we have
\begin{equation}
\|\rho_j(g)-\rho_j(h)\|_{2,D_j}
\ge \frac{s}{2}
\qquad
(g,h\in G_j,\ g\neq h).
\label{eq:rho-separation}
\end{equation}
By Remark~\ref{rem:generalpredicate}, the conclusion also holds when \(\varphi\) is an ultraproduct predicate arising from uniformly equicontinuous coordinate predicates, provided it is \(L\)-almost near.
\end{cor}

\begin{proof}
By Theorem~\ref{thm:pseudofinite}, there are finite
groups $G_j$, numbers $\varepsilon_j>0$ with
$\lim_{\mathcal U}\varepsilon_j=0$, and
$\varepsilon_j$-almost representations $f_j:G_j\to\mathrm U(n_j)$ 
such that $G=\prod_{\mathcal U}f_j(G_j)$. 
Moreover, by \eqref{eq:fjestimate}, there are
numbers $e_j\ge0$ with
$\lim_{\mathcal U}e_j=0$ such that
\[
\|f_j(g)-f_j(h)\|_{2,n_j}\ge s-e_j
\qquad(g,h\in G_j,\ g\ne h)
\]
for $\mathcal U$-almost all $j$. Choose $I\in\mathcal U$ such that for every $j\in I$ the preceding
separation bound holds and $\varepsilon_j$ is small enough for
Fact~\ref{fact:flexible}. For $j\notin I$, replace $G_j$ by the
trivial group, take $f_j$ to be the trivial representation into
$\mathrm U(n_j)$, and set $\varepsilon_j$ equal to a fixed
sufficiently small positive number. Since $I\in\mathcal U$, we still
have $G=\prod_{\mathcal U}f_j(G_j)$ and
$\lim_{\mathcal U}\varepsilon_j=0$. The separation bound is vacuous
for $j\notin I$.

Applying Fact~\ref{fact:flexible}, we obtain
integers $D_j\ge n_j$, isometries
$J_j:\mathbb C^{n_j}\to\mathbb C^{D_j}$,
and unitary representations
$\rho_j:G_j\to\mathrm U(D_j)$ such that 
$\frac{D_j}{n_j}\le1+a(\varepsilon_j)$ 
and $\sup_{g\in G_j}
\|\rho_j(g)-\widehat f_j(g)\|_{2,D_j}
\le b(\varepsilon_j)$, 
where $\widehat f_j(g)
=J_jf_j(g)J_j^*+(1-J_jJ_j^*)$. 
Since $a(\varepsilon),b(\varepsilon)\to0$
as $\varepsilon\to0$, we have
\[
\lim_{\mathcal U}\frac{D_j}{n_j}=1,
\qquad
\lim_{\mathcal U}b(\varepsilon_j)=0.
\]

Set $p_j=J_jJ_j^*$. Since
$\operatorname{tr}_{D_j}(p_j)
=\frac{n_j}{D_j}\longrightarrow_{\mathcal U}1$, the maps $x_j\mapsto J_jx_jJ_j^*$ induce
the inverse corner isomorphism
$\Phi:\M_{\mathbf n}\to\M_{\mathbf D}$. 

For every $g\in G_j$,
\[
\|\rho_j(g)-J_jf_j(g)J_j^*\|_{2,D_j}
\le b(\varepsilon_j)
+\sqrt{1-\frac{n_j}{D_j}}.
\]
The right-hand side converges to zero along
$\mathcal U$, uniformly in $g$.
It follows that
\[
\Phi(G)=\prod_{\mathcal U}\rho_j(G_j).
\]

Finally, for distinct $g,h\in G_j$,
\[
\begin{aligned}
&\|\rho_j(g)-\rho_j(h)\|_{2,D_j}\\
&\quad\ge
\|\widehat f_j(g)-\widehat f_j(h)\|_{2,D_j}
-2b(\varepsilon_j)\\
&\quad=
\sqrt{\frac{n_j}{D_j}}\,
\|f_j(g)-f_j(h)\|_{2,n_j}
-2b(\varepsilon_j)\\
&\quad\ge
\sqrt{\frac{n_j}{D_j}}\,(s-e_j)
-2b(\varepsilon_j).
\end{aligned}
\]
The last expression converges to $s$
along $\mathcal U$, independently of
$g,h$. Consequently, for
$\mathcal U$-almost all $j$,
\[
\|\rho_j(g)-\rho_j(h)\|_{2,D_j}
\ge\frac{s}{2}
\qquad(g,h\in G_j,\ g\ne h).
\]
This proves \eqref{eq:rho-separation}.
\end{proof}

\begin{remark}\label{rem:tracezero}
We note that the hypothesis of Theorem \ref{thm:pseudofinite} and Corollary \ref{cor:pseudofinite} on a uniform lower bound on distances between distinct elements of $G$ is implied by the assumption that each nonidentity element of $G$ has trace $0$.  
\end{remark}

We are nearly ready to prove Theorem \ref{thm:frame}.  First, we need one simple lemma:

\begin{lem}\label{rec:dominant}
For each $j$, suppose that $H_j$ is a finite group and $\rho_j:H_j\to \mathrm{U}(D_j)$ is a unitary representation of $H_j$.  Set $C_j:=\rho_j(H_j)'\cap M_{D_j}(\C)$.  Suppose that $\prod_\cU C_j=\mathbb C$.  Then for $\cU$-almost all $j$, $\rho_j$ has a unique irreducible summand whose dimension is at least $D_j/2$ (and whose multiplicity is necessarily $1$).  Moreover, the support $p_j$ of this summand is a central projection in $C_j$ satisfying $\lim_\cU \tr_{D_j}(p_j)=1$. 
\end{lem}

\begin{proof}
Set
\[
 \mu_j:=\sup
       \{\min(\tr_{D_j}(p),1-\tr_{D_j}(p)) \ : \ p\in C_j \text{ a projection}\}.
\]
Since $\prod_\cU C_j=\bb C$, we have that $\lim_\cU \mu_j=0$.

Decompose $\rho_j$ into irreducibles and decompose its
identity as a sum of mutually orthogonal minimal projections of
$C_j$.  These projections are the supports of individual
irreducible summands.  Every sum of such projections is a
projection of $C_j$, so its trace is either at most $\mu_j$
or at least $1-\mu_j$.  Suppose $\mu_j<1/3$ (which happens for $\cU$-almost all $j$) and yet, towards a contradiction, that none of the aforementioned minimal projections has trace at least $1-\mu_j$.  Then every
such minimal projection has trace at most $\mu_j$.  Adding traces until
their sum first exceeds $\mu_j$ produces a projection in $C_j$ whose trace is in
$(\mu_j,2\mu_j]$, which is thus strictly between $\mu_j$ and
$1-\mu_j$, a contradiction.  Thus, one of the minimal projections has trace at least $1-\mu_j$.
It is unique since $\mu_j<1/2$.  For the same reason, the corresponding irreducible summand has multiplicity $1$, which also establishes that its support is central in $C_j$.
\end{proof}

We are now ready for the proof of Theorem \ref{thm:frame}:

\begin{proof}[Proof of Theorem \ref{thm:frame}]
By Theorem \ref{thm:pseudofinite} and Remark \ref{rem:tracezero}, there are finite groups $G_j$
and $\epsilon_j$-almost representations $f_j:G_j\to\mathrm U(n_j)$ with $\lim_\cU \epsilon_j=0$ for which $G=\prod_\cU f_j(G_j)$.

Write $F:=(F_j)_\cU$, with each $F_j\in \mathrm{U}(n_j)$. Since $\lim_{\mathcal{U}}\|F_j^2-1\|_{2,n_j}=0$, replacing $F_j$ by $F_j':=2\cdot 1_{[0,\infty)}(\frac{1}{2}(F_j+F_j^*))-1$, which satisfies $\lim_{\cU}\|F_j-F_j'\|_{2,n_j}=0$, we may assume that  $F_j^2=1$.  For $j\in \N$ and $g\in G_j$, set $A_{j,g}:=f_j(g)$ and $B_{j,g}:=F_jf_j(g)F_j$.
Set
\[
 \gamma_j:=\sup_{g,h\in G_j}
       \|[A_{j,g},B_{j,h}]\|_2.
\]
The fact that $[G,FGF^*]=0$ implies that $\lim_\cU \gamma_j=0$.  

Form the semidirect product
\[
 H_j:=(G_j\times G_j)\rtimes C_2,
\]
where the nonidentity element of $C_2$ swaps the two copies of $G_j$.
Next define $f'_j:H_j\to \mathrm{U}(n_j)$ by
\[
 f'_j(g,h,t)=A_{j,g}B_{j,h}F_j^{\,t},
 \qquad t\in\{0,1\}.
\]

It is straightforward to check that $f'_j$ is a $(2\epsilon_j+\gamma_j)$-
almost representation. 

By Fact \ref{fact:flexible}, we obtain integers $D_j\geq n_j$ satisfying $\mathbf{D}/\mathbf{n}=1$ and representations
\[
 \rho_j:H_j\to\mathrm{U}(D_j)
\]

such that $\lim_\cU \sup_{g\in H_j}\|\rho_j(g)-\hat f'_j(g)\|_2=0$.

Since $\mathbf D/\mathbf n=1$, the corner
isomorphism introduced in Subsection \ref{subsec:proofofmaintheorem}
yields a trace-preserving $*$-isomorphism
$\Phi_1:\M_{\mathbf n}\to\M_{\mathbf D}$. We choose $\Phi_1$ using the isometries
$J_j:\C^{n_j}\to\C^{D_j}$ provided by
Fact~\ref{fact:flexible}; explicitly, 
\[\Phi_1((u_j)_{\cU})=(J_ju_jJ_j^*)_{\cU},\qquad (u_j)_{\cU}\in \M_{\mathbf{n}}.\]

We then have the following equations:
\begin{align*}
 \Phi_1(G)&=\prod_\cU \rho_j(G_j\times\{e\}\times \{0\}),\\
 \Phi_1(FGF^*)&=\prod_\cU \rho_j(\{e\}\times G_j\times \{0\}),\\
 \Phi_1(F)&=(\rho_j(e,e,1))_\cU.
\end{align*}
View $G_j\times G_j$ as a subgroup of $H_j$ by identifying it with $G_j\times G_j\times \{0\}$.  
Let
\[
 C_j=\rho_j(G_j\times G_j)'
          \cap M_{D_j}(\C).
\]
The fact that $(G\cup FGF^*)'\cap \M_{\mathbf{n}}=\bb C$ implies that
$\prod_{\cU}C_j=\C1$.  Write $\rho_j\downarrow (G_j\times G_j)$ for the restriction of $\rho_j$ to $G_j\times G_j$.  Lemma~\ref{rec:dominant}
therefore gives a central projection $p_j\in C_j$ such that $\lim_\cU \tr_{D_j}(p_j)=1$ and such that $(\rho_j\downarrow (G_j\times G_j))|p_j(\bb C^{D_j})$ is irreducible and occurs with multiplicity one.  It follows that $$(\rho_j\downarrow (G_j\times G_j))|p_j(\bb C^{D_j})\cong \pi_j\otimes\sigma_j,$$
with $\pi_j$ and $\sigma_j$ irreducible representations
of $G_j$ on spaces $V_j$ and $W_j$ respectively.

Write $T_j=\rho_j(e,e,1)$, so $\Phi_1(F)=(T_j)_\cU$.  Notice that $T_j\rho_j(g,h,0)T_j^*=\rho_j(h,g,0)$ for all $g,h\in G_j$.  Since conjugation by $T_j$ sends $p_j$ to a central support of another irreducible
summand of the same dimension, uniqueness of the summand of
weight greater than $1/2$ implies that $T_jp_jT_j^*=p_j$.
Consequently, $T_j|p_j(\bb C^{D_j})$ implements an equivalence
\[
 \pi_j\otimes\sigma_j\cong\sigma_j\otimes\pi_j.
\]
We next claim that $\pi_j$ is equivalent to $\sigma_j$.  To see this, note that, taking characters we have $\chi_{\pi_j}(g)\chi_{\sigma_j}(h)=\chi_{\sigma_j}(g)\chi_{\pi_j}(h)$ for all $g,h\in G_j$.  Applying this with $h=e$ gives $\chi_{\pi_j}(g)\dim(W_j)=\chi_{\sigma_j}(g)\dim (V_j)$.  Letting $c:=\dim(V_j)/\dim(W_j)$, we see that $\chi_{\pi_j}=c\cdot \chi_{\sigma_j}$.  Since $\pi_j$ and $\sigma_j$ are both irreducible, we have that $\langle \pi_j,\pi_j\rangle_{G_j}=\langle \sigma_j,\sigma_j\rangle_{G_j}=1$, which forces $c=1$ and so $\chi_{\pi_j}=\chi_{\sigma_j}$ and thus $\pi_j$ is equivalent to $\sigma_j$, as claimed.  

Write $d_j=\dim\pi_j$, so $\operatorname{rank}(p_j)=d_j^2$.
After a unitary identification of the two representation spaces,
we may assume that
\[
 \rho_j(g,h,0)|{p_j(\bb C^{D_j})}=\pi_j(g)\otimes\pi_j(h).
\]
Note also that
\[
 \lim_\cU\frac{d_j^2}{n_j}
 =\lim_\cU\frac{\operatorname{rank}(p_j)}{D_j}
       \frac{D_j}{n_j}=1.
\]

Irreducibility of $\pi_j$ gives that the relative commutant and relative bicommutant of $\pi_j(G_j)\otimes 1$ in $M_{d_j}(\C)\otimes M_{d_j}(\C)$ are $1\otimes M_{d_j}(\C)$ and $M_{d_j}(\C)\otimes1$
respectively.

For $v\otimes w\in \bb C^{d_j}\otimes \bb C^{d_j}$, we have
$$(T_j\rho_j(g,h,0)T_j^*)(v\otimes w)=(\pi_j(h)\otimes \pi_j(g))(v\otimes w)=\pi_j(h)v\otimes \pi_j(g)w$$ while
\begin{align*}
 (\operatorname{Flip}_{d_j}\rho_j(g,h,0)\operatorname{Flip}_{d_j}^*)(v\otimes w)&=\operatorname{Flip}_{d_j}\rho_j(g,h,0)(w\otimes v)   \\
 &=\operatorname{Flip}_{d_j}(\pi_j(g)w\otimes \pi_j(h)v)\\
 &=\pi_j(h)v\otimes \pi_j(g)w.
\end{align*}
It follows that 
$(T_j|_{p_j(\bb C^{D_j})})\operatorname{Flip}_{d_j}^{-1}$ 
commutes with
the irreducible representation $(\rho_j\downarrow (G_j\times G_j))|p_j(\bb C^{D_j})$, and so Schur's lemma gives
$T_j|{p_j(\bb C^{D_j})}=\lambda_j\operatorname{Flip}_{d_j}$ for a scalar $\lambda_j$ of modulus one.
Since both $T_j|p_j(\bb C^{D_j})$ and $\operatorname{Flip}_{d_j}$ are involutions, we have that $\lambda_j^2=1$.  It follows that $\operatorname{Ad}(T_j|p_j(\bb C^{D_j}))=\operatorname{Ad}(\operatorname{Flip}_{d_j})$.

Since $\lim_\cU \tr_{D_j}(p_j)=1$, we have the corner isomorphism $\Phi_2:\M_{\mathbf{D}}\to \prod_\cU p_jM_{D_j}(\bb C)p_j$, and since $\operatorname{rank}(p_j)=d_j^2$, we have the isomorphism $\Phi_3:\prod_{\cU}p_jM_{D_j}(\bb C)p_j\to \M_{\mathbf{d}^2}$.  Using Lemmas \ref{rec:internal-algebra} and \ref{rec:internal-group}, we see that the isomorphism $\Phi:\M_{\mathbf{n}}\to \M_{\mathbf{d}^2}$ given by $\Phi:=\Phi_3\circ\Phi_2\circ \Phi_1$ satisfies the conclusion of the theorem.
\end{proof}

\section{Proof of Theorem \ref{thm:affine}}\label{sec:affine}

In this section, we prove Theorem \ref{thm:affine}.  Throughout, we fix a prime power $q\geq 3$ and set
\[
 H_q:=\{(b,a):a\in\mathbb F_q^\times,\ b\in\mathbb F_q\},
 \qquad (b,a)\cdot (e,c):=(b+ae,ac).
\]

In other words, $H_q$ is the group of affine transformations of $\mathbb{F}_q$ by identifying ${(b,a)}$ with the affine transformation $x\mapsto ax+b$.
Throughout this section write
\[
 d=q-1,\qquad n=|H_q|=q(q-1),\qquad
 \mathcal H_q=\ell^2(H_q),\qquad
 \tr_n=n^{-1}\operatorname{Tr}.
\]

The representation theory for $H_q$ is well-known (see \cite[Chapter 25]{james2001representations}). We recall the short argument. Let $\vartheta_q$ be the representation of $H_q$ on the orthogonal
complement of the constants in $\ell^2(\mathbb F_q)$ under the
natural permutation action. It has degree $d=q-1$ and is realized
by real orthogonal matrices. Since the affine action is doubly
transitive, it has two orbits on $\mathbb F_q\times\mathbb F_q$,
so the permutation character has squared norm $2$. The trivial
representation occurs once, so the character of $\vartheta_q$
has squared norm $1$. Thus $\vartheta_q$ is irreducible.

The quotient map $(b,a)\mapsto a$ gives the $d$ distinct linear
characters
\[
 \chi(b,a)=\chi(a),\qquad \chi\in\widehat{\mathbb F_q^\times}.
\]
Together with $\vartheta_q$, they exhaust the irreducible complex
representations of $H_q$, since $d+d^2=q(q-1)=|H_q|$.

Define commuting unitary representations on $\mathcal H_q$ by
\[
 \lambda_q(g)(\delta_x)=\delta_{gx},\qquad
 \rho_q(g)(\delta_x)=\delta_{xg^{-1}}.
\]
Note that $\lambda_q$ is just the left regular representation of $H_q$ and $\rho_q$ is the right regular representation of $H_q$.  We will use the same symbols $\lambda_q$ and $\rho_q$ to denote the linear
extensions of these maps to $\bb C[H_q]$.  

We note that the decomposition of $\lambda_q$ into irreducible representations is given by 
\[
\lambda_q\cong \bigoplus_{\chi\in\widehat{\mathbb F_q^\times}}\chi
                 \ \oplus\ d\vartheta_q.
\]
Define also the inversion unitary $F_q$ by
\[
 F_q(\delta_x)=\delta_{x^{-1}}.
\]
The following equations are then readily verified:
\begin{equation}
\label{aff:flip-identities}
 F_q^*=F_q,\qquad F_q^2=1,\qquad
 F_q\lambda_q(g)F_q=\rho_q(g),\qquad [\lambda_q(g),\rho_q(h)]=0.
\end{equation}

Set
\[
 A_q:=\lambda_q(\C[H_q]),\qquad B_q:=\rho_q(\C[H_q]). 
\]
It is standard fare that $A_q'\cap M_n(\C)=B_q$ and $B_q'\cap M_n(\C)=A_q$, whence we have $(A_q\cup B_q)'\cap M_n(\bb C)=Z(A_q)$.
We recall the following standard piece of representation
theory. Recall that $\widehat{H_q}$ denotes the set of equivalence classes of irreducible representations of $H_q$.  We will often identify an element $\pi\in\widehat{H_q}$ with a corresponding representative
$\pi:H_q\to\mathrm U(V_\pi)$ of the class.  By linear extension, we obtain a unital
$*$-homomorphism $\pi:\mathbb C[H_q]\to\operatorname{End}(V_\pi)$. 
Taking the direct sum of these maps gives a
vector space isomorphism
\[
\mathbb C[H_q]
\cong
\bigoplus_{\pi\in\widehat{H_q}}
\operatorname{End}(V_\pi).
\]
This isomorphism is in fact an isomorphism
of $\mathbb C[H_q\times H_q]$-modules,
where $H_q\times H_q$ acts on
$\mathbb C[H_q]$ via
\[
(g,h)\longmapsto\lambda_q(g)\rho_q(h),
\]
and acts on $\operatorname{End}(V_\pi)$ by
\[
X\longmapsto\pi(g)X\pi(h^{-1}).
\]
Moreover, identifying
$\operatorname{End}(V_\pi)$ with
$V_\pi\otimes V_\pi^*$, this action becomes
the outer tensor product representation
$\pi\boxtimes\pi^\vee$, where
$\pi^\vee$ denotes the dual
representation of $\pi$ on $V_\pi^*$.
Let $p_q\in M_n(\C)$ be the projection which corresponds to the projection onto the summand $\operatorname{End}(V_{\vartheta_q})$.  Note that $\tr_{n}(p_q)=d^2/n=1-1/q$ and that $p_q\in A_q\cap B_q$ is a central projection.  Moreover, since $\vartheta_q$ is a real, orthogonal representation, we have that $\vartheta_q\cong \vartheta_q^{\vee}$.  It follows that there is a unitary map $w_q\colon p_q\mathcal{H}_q\to \C^d\otimes \C^d$ which intertwines $\lambda_q|p_q\cal H_q$ and $\vartheta_q\otimes 1$ as well as $\rho_q|p_q\cal H_q$ and $1\otimes \vartheta_q$.  Moreover, $w_q$ intertwines $F_q|p_q\cal H_q$ and $\operatorname{Flip}_d$.  In particular, $p_qA_q$ gets identified with $M_d(\C)\otimes 1$ and $p_qB_q$ gets identified with $1\otimes M_d(\C)$.

We can now begin the proof of Theorem \ref{thm:affine}.  The definable subgroup $G$ of $\mathrm{U}(\M_\mathbf{n})$ will be the pseudofinite group $\prod_\cU \lambda_{q_j}(H_{q_j})$.  To show that $G$ is defined by the $L$-almost near formula $\varphi_{\operatorname{aff}}$, we proceed in two steps.  The first step is to bound the distance of an element of $\mathrm{U}(A_q)$ to $\lambda_q(H_q)$.

Consider the unital $*$-homomorphism $\mathcal D:\bb C[H_q]\to \C[H_q]\otimes \C[H_q]$ given by $$\mathcal D\left(\sum_{g\in H_q} c_g g\right):=\sum_{g\in H_q}c_g(g\otimes g).$$

If $a\in \mathrm{U}(\bb C[H_q])$, then both
$\mathcal D(a)$ and $a\otimes a$ are unitary.  Let $\tau_{\reg}$ denote the usual normalized trace on $\bb C[H_q]$, that is, $\tau_{\reg}(e)=1$ and $\tau_{\reg}(g)=0$ if $g\in H_q\setminus \{e\}$.  In other words, $\tau_{\reg}(x)=\tr_n(\lambda_q(x))$ for all $x\in \bb C[H_q]$.  Set $\omega_q:=\tau_{\reg}\otimes \tau_{\reg}$.  It follows that
\begin{align*}
 \omega_q(\mathcal D(a)^*(a\otimes a))
 &=\sum_{g,h,k}\overline{c_g}c_hc_k
     \tau_{\mathrm{reg}}(g^{-1}h)\tau_{\mathrm{reg}}(g^{-1}k)\\
 &=\sum_g|c_g|^2c_g.
\end{align*}
Consequently,
\begin{equation}
\label{aff:exact-fourier-identity}
 \|\mathcal D(a)-a\otimes a\|_{2,\omega_q}^2
 =2-2\operatorname{Re}\sum_g|c_g|^2c_g.
\end{equation}
On the other hand,
\[
 \|\lambda_q(a)-\lambda_q(h)\|_{2,n}^2=2-2\operatorname{Re}c_h.
\]
The real part of the weighted sum in
\eqref{aff:exact-fourier-identity} is at most
$\max_h\operatorname{Re}c_h$, since the nonnegative weights
$|c_g|^2$ sum to one.  It follows that
\[
 d(\lambda_q(a),\lambda_q(H_q))
 \le\|\mathcal D(a)-a\otimes a\|_{2,\omega_q}.
\]

The term appearing on the right is easier to bound if we switch to a different state on $\bb C[H_q]\otimes \bb C[H_q]$.  Indeed, let $\Omega_q:H_q\times H_q\to \mathrm{U}(n)$ be the representation defined by $\Omega_q(g,h)=\lambda_q(g)\rho_q(h)$, which is indeed a representation since the left and right regular representations commute.  Extend $\Omega_q$ to $\bb C[H_q\times H_q]$ by linearity.  We get an associated state $\sigma_q$ on $\bb C[H_q]\otimes \bb C[H_q]$ given by $\sigma_q(a\otimes b)=\tr_{{n}}(\lambda_q(a)\rho_q(b))$.  Note then that, for $a\in \bb C[H_q]$, we have $\Omega_q(\cal D(a)-a\otimes a)=\Delta_q(a)-\lambda_q(a)\rho_q(a)$, where $\Delta_q:H_q\to \mathrm{U}(n)$ is the representation of $H_q$ given by $\Delta_q(g)=\Omega_q(g,g)=\lambda_q(g)\rho_q(g)$.  It follows that $$\|\cal D(a)-a\otimes a\|_{2,\sigma_q}=\|\Delta_q(a)-\lambda_q(a)\rho_q(a)\|_2.$$  We now show two things:  (1) the norm $\|\Delta_q(a)-\lambda_q(a)\rho_q(a)\|_2$ can be easily approximated from above and, (2) the norms calculated with respect to the two states $\omega_q$ and $\sigma_q$ are not all that different.

To explain the first item, we need to describe the decomposition of $\Delta_q$ into irreducible representations.  Perhaps this is well-known, but we give a proof for completeness:

\begin{lem}
The decomposition of $\Delta_q$ into irreducible representations is given by
\[
 \Delta_q\cong q\,1
       \ \oplus\!\bigoplus_{\substack{\chi\in\widehat{\mathbb F_q^\times}\\
                                      \chi\ne1}}\chi
       \ \oplus\ (d-1)\vartheta_q.
\]
Consequently, $\lambda_q$ and $\Delta_q$ have equivalent invariant
subspaces of dimension $d^2$ (and thus of codimension $d$).
\end{lem}

\begin{proof}
For $g\in H_q$, we have that $\Delta_q(g)(x)=gxg^{-1}$ for all $x\in H_q$.  If $\chi_{\Delta_q}$ denotes its character, then $\chi_{\Delta_q}(g)=|C_{H_q}(g)|$.  This character thus has value $n$ at the
identity, $q$ on nonidentity translations, and $d$ on the remaining
elements.  Its inner product with a linear character $\chi$ is
\[
 \frac1{qd}\left(qd+dq+qd\sum_{a\ne1}\overline{\chi(a)}\right)
 =2+\sum_{a\ne1}\overline{\chi(a)}.
\]
This is $2+(q-2)=q$ for the trivial character, and $2-1=1$ for each
nontrivial linear character.  Its inner product with $\chi_{\vartheta_q}$ is
\[
 \frac1{qd}\left(qd\cdot d+dq\cdot(-1)\right)=d-1.
\]
The established decomposition follows from these calculations.
\end{proof}

We can now achieve the first item:

\begin{lem}
\label{aff:intertwiner}
There is a unitary $W_q\in B(\mathcal H_q)$ such that, for all $a\in \mathrm{U}(\bb C[H_q])$, we have
\begin{equation}
\label{aff:whole-algebra-intertwining}
 \|W_q\lambda_q(a)W_q^*-\Delta_q(a)\|_2
 \le\frac2{\sqrt q}.
\end{equation}
\end{lem}

\begin{proof}
Since $\lambda_q$ and $\Delta_q$ have equivalent invariant subspaces of dimension $d^2$, we may choose a unitary intertwiner
between them and extend it to a unitary $W_q$ using any unitary
between their orthogonal complements (whose common dimension is $d$).
On the image of the common invariant subspace,
$W_q\lambda_q(a)W_q^*$ and $\Delta_q(a)$ agree for every $a \in \bb C[H_q]$.  If, moreover, $a$ is unitary, both operators are
unitary, so their difference has operator norm at most $2$ and is
supported on a subspace of dimension $d$.  It follows that
$$\|W_q\lambda_q(a)W_q^*-\Delta_q(a)\|^2_2
 \le 4d/n=4/q,$$
as claimed.
\end{proof}

We now prove the second item:

\begin{lem}
\label{aff:state-comparison}
As linear functionals on the finite-dimensional $C^*$-algebra
$\C[H_q]\otimes\C[H_q]$, we have
\begin{equation}
\label{aff:state-bound}
 \|\sigma_q-\omega_q\|\le \frac4q.
\end{equation}
Consequently, for every element $z$ of that algebra with
$\|z\|\le2$, we have
\begin{equation}
\label{aff:norm-transfer}
 \omega_q(z^*z)^{1/2}
 \le \sigma_q(z^*z)^{1/2}+\frac4{\sqrt q}.
\end{equation}
\end{lem}

\begin{proof}
Set $P:=d^2/n$, which is the weight with respect to $\tau_{\reg}$ of the $\vartheta_q$ block. It follows that the weights of $\omega_q$ and $\sigma_q$ on the $(\vartheta_q,\vartheta_q)$ block are $P^2$ and $P$ respectively.  Let $\nu$ be the normalized trace of the $(\vartheta_q,\vartheta_q)$ block (which is the same for both tracial states), extended to $0$ on the other blocks.  It follows that $\omega_q$ restricted to the $(\vartheta_q,\vartheta_q)$ block is $P^2\nu$ while $\sigma_q$ restricted to the $(\vartheta_q,\vartheta_q)$ block is $P\nu$.  Next note that $\omega_q-P^2\nu$ is positive since it is simply $\omega_q$ restricted to the other blocks; moreover, $\|\omega_q-P^2\nu\|=(\omega_q-P^2\nu)(1)=1-P^2$.  Similarly, $\sigma_q-P^2\nu$ is positive as it coincides with $(P-P^2)\nu$ plus $\sigma_q$ on the other blocks; moreover, $\|\sigma_q-P^2\nu\|=(\sigma_q-P^2\nu)(1)=1-P^2$.   It follows that $$\|\omega_q-\sigma_q\|\leq 2(1-P^2)=2(2/q-1/q^2)\leq 4/q,$$
yielding (\ref{aff:state-bound}).

If $\|z\|\le2$, then $\|z^*z\|\le4$, and therefore
\[
 \omega_q(z^*z)\le\sigma_q(z^*z)+\frac{16}{q}.
\]
Taking square roots and using $\sqrt{s+t}\le\sqrt s+\sqrt t$
proves \eqref{aff:norm-transfer}.
\end{proof}

We can now summarize the above discussion and prove a weak definability result, estimating how far away elements of $\mathrm{U}(A_q)$ are from $\lambda_q(H_q)$.  Set $$\gamma_q(u,W_q,F_q):=\|W_quW_q^*-uF_quF_q\|_2.$$

\begin{lem}
\label{aff:fourier-rounding}
For every $a\in\mathrm{U}(\C[H_q])$, we have
\begin{equation}
\label{aff:algebra-rounding}
 d(\lambda_q(a),\lambda_q(H_q))
 \le
 \gamma_q(\lambda_q(a),W_q,F_q)
 +\frac6{\sqrt q}.
\end{equation}
\end{lem}

\begin{proof}
Fix $a\in \mathrm{U}(\bb C[H_q])$.  As established above, we have
$$d(\lambda_q(a),\lambda_q(H_q))\leq \|\cal D(a)-a\otimes a\|_{2,\omega_q}.$$
Using  Lemmas~\ref{aff:intertwiner} and \ref{aff:state-comparison} together with the fact already mentioned above, namely that $\Omega_q(\cal D(a)-a\otimes a)=\Delta_q(a)-\lambda_q(a)\rho_q(a)$, we have
\begin{align*}
 d(\lambda_q(a),\lambda_q(H_q))
 &\le\|\Delta_q(a)-\lambda_q(a)\rho_q(a)\|_2+4q^{-1/2}\\
 &\le\|W_q\lambda_q(a)W_q^*-\lambda_q(a)\rho_q(a)\|_2
           +6q^{-1/2}.
\end{align*}
Since $\rho_q(a)=F_q\lambda_q(a)F_q$, the lemma is proven.
\end{proof}

We now must extend the previous lemma to include all unitaries in $\mathrm{U}(n)$.  To do this, we recall the following fact of Ioana (see \cite[Proposition 4.1]{ioana2025almost}).

\begin{fact}
\label{aff:external-expander}
There is a constant $\eta>0$ with the following property:
For every pair of positive integers $k,m$, there are two unitaries
$Z_1,Z_2\in\mathrm U(M_k(\C)\otimes1_m)$ such that
\[
 \|x-E_{1_k\otimes M_m(\C)}(x)\|_2
 \le\eta\bigl(\|[Z_1,x]\|_2+\|[Z_2,x]\|_2\bigr)
 \qquad(x\in M_k(\C)\otimes M_m(\C)).
\]
\end{fact}

Applying this result, with the tensor factors interchanged, to the ``dominant'' 
$d\times d$ block, we may take  $v_1,v_2\in\mathrm U(d)$ such that
\begin{equation}
\label{aff:selector-gap}
 \|x-E_{M_d(\C)\otimes1}(x)\|_{2,p_q}
 \le\eta\sum_{i=1}^2\|[x,1\otimes v_i]\|_{2,p_q}
\end{equation}
holds for all $x\in M_d(\C)\otimes M_d(\C)$.
Here $\|\cdot\|_{2,p_q}$ denotes the normalized norm on the corner.
Define $b_{i,q}\in \mathrm{U}(B_q)$ to be $1\otimes v_i$ on $p_q\cal H_q$ and the
identity on each complementary one-dimensional block.
For an arbitrary unitary $u\in\mathrm U(n)$, set
\[
\beta_q(u,b_{1,q},b_{2,q})=\max_{i=1,2}\|[u,b_{i,q}]\|_2.
\]

We can now measure the distance from a unitary in $\mathrm{U}(n)$ to one in $\mathrm{U}(A_q)$:

\begin{lem}
\label{aff:selector-rounding}
For every $u\in\mathrm U(n)$, there is $v\in\mathrm U(A_q)$ such that
\begin{equation}
\label{aff:selector-unitary-bound}
 \|u-v\|_2\le2\sqrt2\eta\,\beta_q(u,b_{1,q},b_{2,q})+\frac2{\sqrt q}.
\end{equation}
\end{lem}

\begin{proof}
Let $p=p_q$ and let $E_{A_q}$ denote the
trace-preserving conditional expectation.
Put
\[
z=E_{A_q}(u),
\qquad
y=E_{A_q}(pup)=pzp,
\]
where the last equality follows from the
$A_q$-bimodularity of $E_{A_q}$.

The restriction of $E_{A_q}$ to $pM_n(\C)p$
is the expectation onto $M_d(\C)\otimes1$.
Multiplying \eqref{aff:selector-gap} by
$\tr_{n}(p)^{1/2}$ gives the same inequality
in the ambient norm.
Since $p$ commutes with each $b_{i,q}$,
we have
\[
[pup,pb_{i,q}]=p[u,b_{i,q}]p.
\]
It follows that
\[
\|pup-y\|_{2,n}
=\|pup-E_{A_q}(pup)\|_{2,n}
\le 2\eta\beta_q(u).
\]

Since
\[
\begin{aligned}
\|pup\|_{2,n}^2
&=\tr_{n}(p)-\|(1-p)up\|_{2,n}^2\\
&\ge (1-q^{-1})-q^{-1},
\end{aligned}
\]
we have
\[
\|u-pup\|_{2,n}^2
=1-\|pup\|_{2,n}^2
\le\frac2q.
\]

Now, since $z=E_{A_q}(u)$ is the
orthogonal projection of $u$ onto $A_q$
in $L^2(M_n(\C),\tr_{{n}})$ and $y\in A_q$,
we have
\[
\|u-z\|_{2,n}^2\le\|u-y\|_{2,n}^2.
\]
Moreover, $u-pup$ is orthogonal to
$pup-y\in pM_n(\C)p$. Consequently,
\[
\begin{aligned}
\|u-z\|_{2,n}^2
&\le\|u-y\|_{2,n}^2\\
&=\|u-pup\|_{2,n}^2
  +\|pup-y\|_{2,n}^2\\
&\le4\eta^2\beta_q(u)^2+\frac2q.
\end{aligned}
\]
Thus
\begin{equation}
\label{aff:expectation-error}
\|u-z\|_{2,n}^2
\le4\eta^2\beta_q(u)^2+\frac2q.
\end{equation}

Since $A_q$ is a finite von Neumann algebra,
the partial isometry in the polar decomposition
of $z$ extends to a unitary $v\in A_q$.
Thus $z=v|z|$.

With this extension,
\[
 \|z-v\|_2^2
 =\tr_{n}((1-|z|)^2)
 \le\tr_{n}(1-|z|^2)=\|u-z\|_2^2.
\]
Since $u-z$ is orthogonal to $z-v\in A_q$, we have
\[
 \|u-v\|_2^2
 =\|u-z\|_2^2+\|z-v\|_2^2
 \le2\|u-z\|_2^2.
\]
Using \eqref{aff:expectation-error} and
$\sqrt{s+t}\le\sqrt s+\sqrt t$ gives
\eqref{aff:selector-unitary-bound}.
\end{proof}

Now set $\varphi_{\operatorname{aff}}(u,b_{1,q},b_{2,q},W_q,F_q)=\max\{\beta_q(u,b_{1,q},b_{2,q}),\gamma_q(u,W_q,F_q)\}$.  The following is our main definability result:

\begin{prop}
\label{aff:ambient-rounding}
Set
\[
 L=1+8\sqrt2\eta,\qquad \varepsilon_q=14q^{-1/2}.
\]
For every $u\in\mathrm U(n)$, we have
\begin{equation}
\label{aff:ambient-distance}
 d(u,\lambda_q(H_q))\le L\varphi_{\operatorname{aff}}(u,b_{1,q},b_{2,q},W_q,F_q)+\varepsilon_q.
\end{equation}
Moreover,
\begin{equation}
\label{aff:group-defect}
 \max_{g\in H_q}\varphi_{\operatorname{aff}}(\lambda_q(g))\le2q^{-1/2}.
\end{equation}
In particular, we have
\begin{equation}
\label{aff:rounding-formula}
 \sup_{u\in\mathrm U(n)}\inf_{v\in\mathrm U(n)}
 \max\bigl(\varphi_{\operatorname{aff}}(v),
       \|u-v\|_2\dminus L\varphi_{\operatorname{aff}}(u)\bigr)
 \le14q^{-1/2}.
\end{equation}
For simplicity, we have suppressed mention of the parameters in the above three displays.
\end{prop}

\begin{proof}
First note that the predicate $\gamma_q$ is $3$-Lipschitz on the unitary group.
Indeed, conjugation by $W_q$ contributes at most $\|u-v\|_2$
and
\[
 \|uF_quF_q-vF_qvF_q\|_2
 \le\|(u-v)F_quF_q\|_2+\|vF_q(u-v)F_q\|_2
 \le2\|u-v\|_2.
\]
Now fix $u\in \mathrm{U}(n)$ and take $v\in\mathrm U(A_q)$ as in
Lemma~\ref{aff:selector-rounding}.  Applying
Lemma~\ref{aff:fourier-rounding} to $v$ gives
\begin{align*}
 d(u,\lambda_q(H_q))
 &\le \|u-v\|_2+\gamma_q(v)+6q^{-1/2}\\
 &\le \gamma_q(u)+4\|u-v\|_2+6q^{-1/2}\\
 &\le \gamma_q(u)+8\sqrt2\eta\,\beta_q(u)+14q^{-1/2}\\
 &\le L\varphi_{\operatorname{aff}}(u)+\varepsilon_q.
\end{align*}
This proves (\ref{aff:ambient-distance}).  Equation \eqref{aff:group-defect} follows from Lemma~\ref{aff:intertwiner} and the fact that $b_{i,q}\in B_q=A_q'$.  To prove \eqref{aff:rounding-formula}, for each $u\in \mathrm{U}(n)$, choose the group element
realizing the finite minimum in \eqref{aff:ambient-distance} and
use it as the witness $v$ in \eqref{aff:rounding-formula}.
The two terms in that maximum are at most $2q^{-1/2}$ and
$14q^{-1/2}$, respectively.
\end{proof}

Now suppose that $\mathbf{q}$ is a sequence of powers of primes, $n_j:=q_j(q_j-1)$, and $d_j:=(q_j-1)$.  For $i=1,2$, set $b_i:=(b_{i,q_j})_\cU$.  Set also $W:=(W_{q_j})_\cU$ and $F:=(F_{q_j})_\cU$, so $b_1,b_2,W,F\in \mathrm{U}(\M_{\mathbf{n}})$.  Proposition \ref{aff:ambient-rounding} immediately implies that the zeroset $G$ of $\varphi_{\operatorname{aff}}(u,b_1,b_2,W,F)$ is $\prod_\cU \lambda_{H_{q_j}}(H_{q_j})$ and that this formula is $L$-almost near.  

We next show that every nonidentity element of $G$ has trace $0$.  For ease of notation, in the rest of this section, set $H_j:=H_{q_j}$.  Note that if $g\in G$ is not the identity, then $g=(\lambda_{H_j}(g_j))_\cU,$ where $g_j$ is not the identity for $\cU$-almost all $j$; for these $j$, $\lambda_{H_j}(g_j)$ has trace $0$, yielding the desired fact.

By \eqref{aff:flip-identities}, for all $g,h\in H_j$,
\[
[\lambda_{H_j}(g),F_{q_j}\lambda_{H_j}(h)F_{q_j}^*]
=[\lambda_{H_j}(g),\rho_{H_j}(h)]=0.
\]
Passing to the ultraproduct gives $[G,FGF^*]=0$.

We next check that $(G\cup FGF^*)'\cap\M_{\mathbf n}=\C1$.
For each $j$, set
\[
A_j:=\lambda_{H_j}(\C[H_j]),
\qquad
B_j:=\rho_{H_j}(\C[H_j]).
\]
Here, $\rho_{H_j}$ denotes the right regular representation of $H_j$.
Since $B_j=A_j'\cap M_{n_j}(\C)$ and
$A_j=B_j'\cap M_{n_j}(\C)$, we have $
(A_j\cup B_j)'\cap M_{n_j}(\C)=Z(A_j).$

Let $p_j\in Z(A_j)$ be the central projection
corresponding to the irreducible representation
$\vartheta_{q_j}$. Recall that
\[
p_jA_j\cong M_{d_j}(\C),
\qquad
\tr_{n_j}(1-p_j)=\frac1{q_j}.
\]
We claim that $\prod_{\mathcal U}Z(A_j)=\C1$. 
Indeed, let 
$z=(z_j)_{\mathcal U}\in\prod_{\mathcal U}Z(A_j)$,
where $\sup_j\|z_j\|\le K$.
Since $p_jA_j$ is a full matrix algebra,
there is a scalar $c_j$, with $|c_j|\le K$,
such that $p_jz_j=c_jp_j$.

Consequently,
\[
\begin{aligned}
\|z_j-c_j1\|_{2,n_j}
&=\|(z_j-c_j1)(1-p_j)\|_{2,n_j}\\
&\le 2K\tr_{n_j}(1-p_j)^{1/2}\\
&=\frac{2K}{\sqrt{q_j}}
\longrightarrow_{\mathcal U}0.
\end{aligned}
\]
Since the bounded sequence $(c_j)_j$ has an
ultralimit $c\in\C$, it follows that
$z=c1$, proving the claim.

Finally, consider the representation
\[
\Omega_j:H_j\times H_j\to\mathrm U(n_j),
\qquad
\Omega_j(g,h)=\lambda_{H_j}(g)\rho_{H_j}(h).
\]
By the identities
$F_{q_j}\lambda_{H_j}(g)F_{q_j}^*=\rho_{H_j}(g)$
at each coordinate, we have
\[
\langle G,FGF^*\rangle
=\prod_{\mathcal U}\Omega_j(H_j\times H_j).
\]
Applying Lemma~\ref{rec:internal-group},
we obtain
\[
\begin{aligned}
&(G\cup FGF^*)'\cap\M_{\mathbf n}\\
&\quad=
\prod_{\mathcal U}
\left(
\Omega_j(H_j\times H_j)'
\cap M_{n_j}(\C)
\right)\\
&\quad=\prod_{\mathcal U}Z(A_j)
=\C1,
\end{aligned}
\]
as required.

We have now verified that $G$ satisfies the hypotheses of Theorem \ref{thm:frame}.  It remains to verify the two stated properties of the corner isomorphism $\Psi:\M_{\mathbf{n}}\to \M_{\mathbf{d}^2}$.  However, these properties follow immediately from the analogous properties in the finite setting on the dominant block.  This finishes the proof of Theorem \ref{thm:affine}.

\section{Proof of Theorem \ref{thm:special}}\label{sec:special}

In this section, we prove Theorem \ref{thm:special}.  We first prove the easier half of the theorem, namely that $\M_{(2^\mathbf{l})^2}$ has a special unitary for any sequence $\mathbf{l}$.  We begin with the following lemma. 

The unitary $V_H$ below is a version of the 
multiplicative unitary (also called the Kac--Takesaki operator). More precisely,
$V_H=W_H^*$ for the fundamental operator
$(W_Hf)(s,t)=f(s,st)$ of
\cite[page 2 and Proposition 2.4.7]{enockschwartz1992kac}.
This construction goes back to Stinespring
\cite{stinespring1959integration}.
The exact characterization of regular representation unitaries
by the corresponding group-like equation follows from
Tatsuuma duality \cite{tatsuuma1967duality}; see
\cite[Section 0]{woronowicz1996multiplicative}
for an explicit operator formulation.
We include a direct proof because we require the quantitative
estimate \eqref{tar:canonical-rounding}, with constant $1$
uniformly over all finite groups $H$.

\begin{lem}\label{tar:canonical-copier}
Let $H$ be any finite group, set $d:=|H|$, and let $\lambda$ be its left
regular representation on $\ell^2(H)\cong \C^d$.  Consider the unitary $V_H$ on $\ell^2(H)\otimes \ell^2(H)\cong \C^d\otimes \C^d$ defined by
\[
 V_H(\delta_y\otimes\delta_z)=\delta_y\otimes\delta_{yz}.
\]
For every $u\in\mathrm U(d)$, we have
\begin{equation}
 \min_{g\in H}\|u-\lambda(g)\|_2
 \leq \|V_H(u\otimes1)V_H^*-u\otimes u\|_2.
 \label{tar:canonical-rounding}
\end{equation}
Moreover, the set of unitaries in $\mathrm U(d)$ which make the right hand side of (\ref{tar:canonical-rounding}) vanish is precisely $\lambda(H)$. 
\end{lem}

\begin{proof}
First, let us set up some notation.  For $u\in\mathrm U(d)$ and $x,y\in H$, let $u_{x,y}$ be the coefficient of
$\delta_x$ in $u\delta_y$, that is, $u\delta_y=\sum_{x\in H}u_{x,y}\delta_x$.  We then have
\[
 (V_H(u\otimes1)V_H^*)(\delta_y\otimes\delta_z)
 =\sum_{x\in H}u_{x,y}(\delta_x\otimes\delta_{xy^{-1}z}).
\]
Consequently, setting $\Delta(u):=V_H(u\otimes1)V_H^*$, we have
\begin{align*}
 \tr_{d^2}\big(\Delta(u)^*(u\otimes u)\big)
 &=\frac1{d^2}\sum_{x,y,z\in H}
       |u_{x,y}|^2 u_{xy^{-1}z,z}\\
 &=\sum_{g\in H}P_g c_g,
\end{align*}
where
\[
 P_g=\frac1d\sum_{y\in H}|u_{gy,y}|^2,\qquad
 c_g=\frac1d\sum_{z\in H}u_{gz,z}
     =\tr_d(\lambda(g)^*u).
\]
Note that the numbers $P_g$ are nonnegative and sum to $1$.  Since both $\Delta(u)$
and $u\otimes u$ are unitary, we have
\begin{align*}
 \|\Delta(u)-u\otimes u\|_2^2
 &=2-2\operatorname{Re}\sum_gP_gc_g\\
 &\geq2-2\max_g\operatorname{Re}c_g
 =\min_g\|u-\lambda(g)\|_2^2.
\end{align*}

This proves \eqref{tar:canonical-rounding}.

For any $g,y,z\in H$, we have
\begin{align*}
(V_H(\lambda(g)\otimes 1)V_H^*)(\delta_y\otimes \delta_z)&=(V_H(\lambda(g)\otimes 1))(\delta_y\otimes \delta_{y^{-1}z})\\
&=V_H(\delta_{gy}\otimes \delta_{y^{-1}z})\\
&=\delta_{gy}\otimes \delta_{gz}\\
&=(\lambda(g)\otimes \lambda(g))(\delta_y\otimes \delta_z),
\end{align*}
proving that $V_H(\lambda(g)\otimes 1)V_H^*=\lambda(g)\otimes \lambda(g)$ and thus $\lambda(H)$ is included in the zeroset of the formula on the right hand side of (\ref{tar:canonical-rounding}); the reverse inclusion follows immediately from (\ref{tar:canonical-rounding}).
\end{proof}

\begin{cor}
\label{tar:canonical-elementary}
Let $H_j=(\mathbb Z/2\mathbb Z)^{l_j}$ and set $d_j:=|H_j|$.  For each $j$, let $V_{H_j}$ be as in Lemma \ref{tar:canonical-copier} and set $V=(V_{H_j})_\cU$.  Then $V$ is a special unitary in $\M_{\mathbf{d}^2}$.  Moreover, we have
\[
 K_{\mathbf{d},V}=\prod_\cU \lambda_{H_j}(H_j)\otimes 1,
\]
and
\[
 C_{\mathbf{d},V}=
 \prod_\cU\lambda_{H_j}(\C[H_j])\otimes 1.
\]
\end{cor}

\begin{proof}
The displayed description of $K_{\mathbf{d},V}$  follows from Lemma~\ref{tar:canonical-copier}
by choosing a nearest element of $\lambda_{H_j}(H_j)$ at each coordinate.  Lemma~\ref{tar:canonical-copier} also yields item (1) of the definition of special unitary while item (2) follows from the fact that each $H_j$ has exponent $2$.  It remains to show that $C_{\mathbf{d},V}$ is abelian, which will follow from the displayed description of $C_{\mathbf{d},V}$ in the statement of the corollary.  Lemma \ref{rec:internal-group} gives that $C_{\mathbf{d},V}=\prod_{\cU}(\lambda_{H_j}(H_j)'\cap M_{d_j}(\C))\otimes 1$.  It remains to note the well-known fact that $\lambda_{H_j}(H_j)'\cap M_{d_j}(\C)=\lambda_{H_j}(\bb C[H_j])$ (using that $H_j$ is abelian).
\end{proof}

We can now finish the ``if'' direction of
Theorem~\ref{thm:special}.
Suppose that $\mathbf d/2^{\mathbf{l}}=1$,
and put $m_j:=2^{l_j}$ and
$k_j:=\min\{d_j,m_j\}$.
Choose a partial isometry
$w_j:\mathbb C^{d_j}\to\mathbb C^{m_j}$
whose initial and final projections both have rank $k_j$.
Since
\[
\lim_{\mathcal U}\frac{k_j}{d_j}
=\lim_{\mathcal U}\frac{k_j}{m_j}=1,
\]
the maps
\[
x_j\longmapsto
(w_j\otimes w_j)x_j(w_j^*\otimes w_j^*)
\]
induce a trace-preserving $*$-isomorphism $\Xi:\mathcal M_{\mathbf d^2}
\longrightarrow\mathcal M_{\mathbf m^2}$. 
Moreover, $\Xi(\mathcal A_{\mathbf d})
=\mathcal A_{\mathbf m}$, and
\[
\Xi\circ\operatorname{Ad}(\operatorname{Flip}_{\mathbf d})
=
\operatorname{Ad}(\operatorname{Flip}_{\mathbf m})
\circ\Xi.
\]
In particular, for $u\in\mathrm U(\mathcal A_{\mathbf d})$,
\[
\Xi(u\otimes u)=\Xi(u)\otimes\Xi(u).
\]
Let $V_{\mathbf m}\in\mathrm U(\mathcal M_{\mathbf m^2})$
be the special unitary given by Corollary \ref{tar:canonical-elementary},
and set $V:=\Xi^{-1}(V_{\mathbf m})$.
Then
\[
\Xi(K_{\mathbf d,V})=K_{\mathbf m,V_{\mathbf m}},
\qquad
\Xi(C_{\mathbf d,V})=C_{\mathbf m,V_{\mathbf m}}.
\]
Since $\Xi$ preserves the $2$-norm, all three
conditions in the definition of a special unitary
are preserved. Thus $V$ is a special unitary
in $\mathcal M_{\mathbf d^2}$, as required.

We now work towards the other direction of Theorem \ref{thm:special}.  We fix a sequence $\mathbf{d}$ and a special unitary $V\in \mathrm{U}(\M_{\mathbf{d}^2})$.  Let $\varphi(u,V,\operatorname{Flip}_\mathbf{d})$ be the predicate on $\mathrm{U}(\A_\mathbf{d})$ given by $\|VuV^*-u\operatorname{Flip}_\mathbf{d}u\operatorname{Flip}_\mathbf{d}\|_2$.  By assumption, this predicate is $1$-almost near.  Moreover,  the zeroset of this predicate is $K_{\mathbf{d},V}$.  Note that any element $u$ of $K_{\mathbf{d},V}$ satisfies $\tr(u)=\tr(u)^2$, whence has trace $0$ if it is not the identity.  Using Remark \ref{rem:generalpredicate} and Corollary \ref{cor:pseudofinite}, we obtain finite groups $H_j$, integers $D_j\geq d_j$ with $\mathbf{D}/\mathbf{d}=1$, and unitary representations $\rho_j:H_j\to \mathrm{U}(D_j)$. 
Let $\Phi_1:\M_{\mathbf d}\to\M_{\mathbf D}$
be the inverse of the canonical corner isomorphism.
Using the tensor squares of the same coordinate
isometries, we obtain an isomorphism
$\Phi_2:\M_{\mathbf d^2}\to\M_{\mathbf D^2}$
satisfying
\[
\Phi_2(u\otimes1)=\Phi_1(u)\otimes1,
\qquad
\Phi_2(u\otimes u)=\Phi_1(u)\otimes\Phi_1(u).
\]
It follows that
$K_{\mathbf D,\Phi_2(V)}
=\Phi_{2}(K_{\mathbf d,V})
=\prod_{\mathcal U}\rho_j(H_j)$ and \[C_{\mathbf{D},\Phi_2(V)}=\prod_\cU (\rho_j(H_j)'\cap M_{D_j}(\C)).\]

We claim that we may further assume that each $H_j$ has
exponent $2$. Indeed, since $K_{\mathbf D,\Phi_2(V)}
=\prod_{\mathcal U}\rho_j(H_j)$ has exponent $2$, we have
\[
\lim_{\mathcal U}
\sup_{g\in H_j}
\|\rho_j(g)^2-1\|_{2,D_j}=0.
\]
Suppose, towards a contradiction, that on a
$\mathcal U$-large set of indices $j$, the group $H_j$
contains an element $g_j$ with $g_j^2\ne e$.
By \eqref{eq:rho-separation} {(where we may take $s=1$)}, we then have
\[
\begin{aligned}
\|\rho_j(g_j)^2-1\|_{2,D_j}
&=\|\rho_j(g_j^2)-\rho_j(e)\|_{2,D_j}\\
&\ge\frac{1}{2},
\end{aligned}
\]
contradicting the preceding limit.
Thus $H_j$ has exponent $2$ for $\mathcal U$-almost
all $j$, and we may assume that this holds for every
$j$. Consequently, each $H_j$ is an elementary
abelian $2$-group, so $H_j\cong(\mathbb Z/2\mathbb Z)^{l_j}$ for some integer $l_j\ge0$. In particular, $|H_j|=2^{l_j}$.  

Our goal is to show that $\mathbf{d}/2^{\mathbf{l}}=1$.  It will suffice to show that $\mathbf{D}/2^{\mathbf{l}}=1$.  Towards that end, for each $j$ and each irreducible character $\chi\in \hat H_j$, we let $m_{j,\chi}$ denote the multiplicity of $\chi$ in $\rho_j$.  Let $S_j:=\{\chi\in \hat H_j \ : \ m_{j,\chi}=1\}$ and let $s_j:=|S_j|$.

\begin{lem}
$\mathbf{s}/\mathbf{D}=1$.
\end{lem}

\begin{proof}
The commutant of $\rho_j(H_j)$ in $M_{D_j}(\C)$ is
$\bigoplus_{\chi:m_{j,\chi}>0}M_{m_{j,\chi}}(\C)$.  We define unitaries $x_j,z_j\in \rho_j(H_j)'\cap M_{D_j}(\C)$ as follows.  If $m_{j,\chi}=m\geq2$, let $X_m,Z_m$
be the direct sums of $\lfloor m/2\rfloor$ copies of the Pauli
unitaries
\[
 X=\begin{pmatrix}0&1\\1&0\end{pmatrix},
 \qquad
 Z=\begin{pmatrix}1&0\\0&-1\end{pmatrix},
\]
and, when $m$ is odd, one additional scalar $1$.  Since $XZ=-ZX$,
\[
 \|[X_m,Z_m]\|_{2,m}^2
 =4\frac{2\lfloor m/2\rfloor}{m}\geq\frac83.
\]
In each multiplicity-one block, choose the scalar $1$.  Taking these
unitaries blockwise gives $x_j,z_j\in\rho_j(H_j)'\cap M_{D_j}(\C)$ with
\[
 \|[x_j,z_j]\|_{2,D_j}^2\geq\frac83 (1-\frac{s_j}{D_j}).
\]
Setting $x:=(x_j)_\cU$ and $z=(z_j)_\cU$, we have $x,z\in C_{{\mathbf{D}},\Phi_2(V)}$, which is abelian.  It follows that
$$0=\|xz-zx\|_2^2=\lim_\cU \|[x_j,z_j]\|_2^2\geq \lim_\cU\frac{8}{3}(1-\frac{s_j}{D_j}).$$  We conclude that $\mathbf{s}/\mathbf{D}=1$.
\end{proof}

We are now left with the task of showing that $\mathbf{s}/2^{\mathbf{l}}=1$.  This is the most delicate part of the proof and will follow from the following statement:

\begin{prop}\label{prop:nearidentity}
$\lim_\cU \sum_{\chi \in \hat {H}_j}|\frac{m_{j,\chi}}{D_j}-\frac{1}{|H_j|}|=0$. 
\end{prop}

To see that Proposition \ref{prop:nearidentity} indeed finishes the proof, it suffices to note that 
$$\left|\frac{s_j}{D_j}-\frac{s_j}{|H_j|}\right|=\sum_{\chi\in S_j}\left|\frac{m_{j,\chi}}{D_j}-\frac{1}{|H_j|}\right|\leq \sum_{\chi \in \hat {H}_j}\left|\frac{m_{j,\chi}}{D_j}-\frac{1}{|H_j|}\right|.$$

We now proceed to prove Proposition \ref{prop:nearidentity}.  To make this more transparent, it helps to recall that the characters of $H_j$ form an orthonormal basis for $\bb C[H_j]$ under the inner product $\langle f_1,f_2\rangle:=\frac{1}{|H_j|}\sum_{g\in H_j}f_1(g)\overline{f_2(g)}$.  So if we consider $\varphi_j\in \bb C[H_j]$ to be the normalized character of $\rho_j$, that is, $\varphi_j(g):=\tr(\rho_j(g))$, then $\varphi_j=\frac{1}{D_j}\sum_{\chi\in \hat{H}_j} m_{j,\chi}\chi$ and so the Fourier coefficients of $\varphi_j$ are $\langle \varphi_j,\chi\rangle=\frac{m_{j,\chi}}{D_j}$ while the Fourier coefficients of $1_{\{e\}}$ (the indicator function of the identity element in $H_j$) are $\langle 1_{\{e\}},\chi\rangle=\frac{1}{|H_j|}\overline{\chi(e)}=\frac{1}{|H_j|}$ since all irreducible characters of an abelian group are one-dimensional.  

Consequently, Proposition~\ref{prop:nearidentity}
asserts that
\[
\lim_{\mathcal U}
\|\widehat{\varphi_j}-\widehat{1_{\{e\}}}\|
_{l^1(\widehat H_j)}=0;
\]
here, for a function $f:H_j\to\mathbb C$, we let
$\widehat f:\widehat H_j\to\mathbb C$ denote
its Fourier transform, given by
\[
\widehat f(\chi)
=\frac1{|H_j|}\sum_{g\in H_j}
f(g)\overline{\chi(g)}.
\]
For each $j$, let $A(H_j)$ denote the Fourier algebra
of $H_j$ (see \cite{Eymard64Fourier} and
\cite[Chapter~2]{Kaniuth-Lau18Fourier}).
Since $H_j$ is finite abelian, $A(H_j)$ consists
of all complex-valued functions on $H_j$, with
pointwise multiplication and norm
\[
\|f\|_{A(H_j)}
=\sum_{\chi\in\widehat H_j}|\widehat f(\chi)|.
\]
It is a unital commutative Banach algebra with
$\|1\|_{A(H_j)}=1$.

Moreover, since
$\widehat{\varphi_j}(\chi)=m_{j,\chi}/D_j$,
we have
\[
\|\varphi_j\|_{A(H_j)}
=\sum_{\chi\in\widehat H_j}
\frac{m_{j,\chi}}{D_j}=1.
\]

By the definition of the Fourier algebra norm,
Proposition~\ref{prop:nearidentity} is thus
equivalent to
\[
\lim_{\mathcal U}
\|\varphi_j-1_{\{e\}}\|_{A(H_j)}=0.
\]
We will prove this by showing that
$\varphi_j$ is an almost idempotent in
$A(H_j)$. A general Banach algebra
perturbation argument will then produce
an idempotent $b_j$ close to $\varphi_j$.
Using the trace-zero property of nonidentity
elements of $K_{\mathbf D,\Phi_2(V)}$,
together with the uniform separation
estimate \eqref{eq:rho-separation}, we will
show that $b_j=1_{\{e\}}$ for
$\mathcal U$-almost all $j$.
The desired result will then follow.

Write $\Phi_{2}(V)=(\hat{V}_j)_\cU$ with each $\hat V_j$ unitary.  Set $$\kappa_j:=\sup_{g\in H_j}\|\hat{V}_j(\rho_j(g)\otimes 1)\hat{V}_j^*-\rho_j(g)\otimes\ \rho_j(g)\|_2,$$ so $\lim_\cU \kappa_j=0$.

The following elementary quantitative estimate
follows from the standard averaging argument
for intertwiners (see, for example,
\cite[Section~5.1]{folland2016course}).
We include the proof to keep track of the
dimension-independent bound.
\begin{lem} \label{lem:approxintertwiner}
Let $H$ be a finite group, let
$\alpha,\beta:H\to\mathrm U(n)$ be unitary
representations, and let $U\in\mathrm U(n)$.
Suppose that
\[
\sup_{g\in H}
\|\beta(g)U-U\alpha(g)\|_{2,n}\le\kappa.
\]
Then:
\begin{enumerate}
\item $\alpha$ and $\beta$ have equivalent
subrepresentations of dimension at least
$(1-\kappa^2/2)n$.

\item If $(a_\pi)_{\pi\in\widehat H}$ and
$(b_\pi)_{\pi\in\widehat H}$ are the irreducible
multiplicities of $\alpha$ and $\beta$,
respectively, then
\[
\frac1n\sum_{\pi\in\widehat H}
(\dim\pi)|a_\pi-b_\pi|\le\kappa^2.
\]
\end{enumerate}
\end{lem}

\begin{proof}
Consider the unitary representation $R$ of $H$ on
the Hilbert space $M_n(\mathbb C)$, equipped
with the normalized Hilbert--Schmidt inner
product, given by
\[
R(g)(X):=\beta(g)X\alpha(g)^*.
\]
The averaging operator $P(X):=\frac1{|H|}\sum_{g\in H}R(g)(X)$ is the orthogonal projection onto the fixed-point
subspace of $R$, which is precisely the space
of intertwiners from $\alpha$ to $\beta$.

Set $T:=P(U)
=\frac1{|H|}\sum_{g\in H}
\beta(g)U\alpha(g)^*$. 
Thus $T\alpha(g)=\beta(g)T$ for every $g\in H$.
Since $P$ is an orthogonal projection and
$\|U\|_{2,n}=1$, we have
\[
\begin{aligned}
\|T-U\|_{2,n}^2
&=1-\|T\|_{2,n}^2\\
&=\frac1{2|H|}\sum_{g\in H}
\|R(g)(U)-U\|_{2,n}^2\\
&=\frac1{2|H|}\sum_{g\in H}
\|\beta(g)U-U\alpha(g)\|_{2,n}^2\\
&\le\frac{\kappa^2}{2}.
\end{aligned}
\]

Let $p$ be the orthogonal projection onto
$\ker T$. Since $Tp=0$ and $U$ is unitary,
\[
\begin{aligned}
\operatorname{tr}_n(p)&=\|Up\|_{2,n}^2
=\|(T-U)p\|_{2,n}^2\\
&\le\|T-U\|_{2,n}^2\le\frac{\kappa^2}{2}.
\end{aligned}
\]
Consequently, $\operatorname{rank}T
\ge(1-\kappa^2/2)n$. 
Write $T=w|T|$ for the polar decomposition
of $T$, and set
\[
p_\alpha=w^*w=s(T^*T),
\qquad
p_\beta=ww^*=s(TT^*).
\]
Here we use $s$ to denote the support projection.  Since $T$ intertwines $\alpha$ and $\beta$,
we have
\[
T^*T\in\alpha(H)',
\qquad
TT^*\in\beta(H)'.
\]
Hence $p_\alpha$ and $p_\beta$ commute with
$\alpha(H)$ and $\beta(H)$, respectively.
In particular, the subspaces
\[
E:=p_\alpha\mathbb C^n,
\qquad
F:=p_\beta\mathbb C^n
\]
are invariant under $\alpha$ and $\beta$.

Moreover, $|T|$ commutes with $\alpha(H)$.
Using $T=w|T|$, we obtain
\[
w\alpha(g)|T|
=w|T|\alpha(g)
=T\alpha(g)
=\beta(g)w|T|.
\]
Since $|T|$ is invertible on $E$, it follows
that
\[
w\alpha(g)|_E=\beta(g)w|_E
\qquad(g\in H).
\]
The restriction $w|_E:E\to F$ is unitary,
and therefore implements an equivalence $\alpha|_E\cong\beta|_F$. 
Since $\dim E=\dim F=\operatorname{rank}T
\ge(1-\kappa^2/2)n$, 
this proves (1).

For (2), the largest possible dimension of
a common subrepresentation of $\alpha$
and $\beta$ is
\[
r:=\sum_{\pi\in\widehat H}
(\dim\pi)\min(a_\pi,b_\pi).
\]
By (1), we have $r\ge(1-\kappa^2/2)n$. 
Since both representations have dimension $n$,
\[
\sum_{\pi\in\widehat H}
(\dim\pi)|a_\pi-b_\pi|
=2n-2r \le\kappa^2 n.
\]
Dividing by $n$ proves (2).
\end{proof}

We now apply Lemma \ref{lem:approxintertwiner} to the representations $\alpha_j,\beta_j:H_j\to \mathrm{U}(D_j^2)$ given by $\alpha_j(g):=\rho_j(g)\otimes 1$ and $\beta_j(g):=\rho_j(g)\otimes \rho_j(g)$.  Noting that each character of $H_j$ is one-dimensional (as $H_j$ is abelian), item (2) of the lemma gives us that
$$\frac{1}{D_j^2}\sum_{\chi \in \hat{H_j}}|D_jm_{j,\chi}-n_{j,\chi}|\leq \kappa_j^2.$$  Here, $n_{j,\chi}$ is the character multiplicity of $\rho_j\otimes \rho_j$. Noting that $\frac{n_{j,\chi}}{D_j^2}=(\hat\varphi_j*\hat\varphi_j)(\chi)$, we have that $\|\hat \varphi_j-\hat \varphi_j*\hat \varphi_j\|_{1}\leq \kappa_j^2$.  Since convolution corresponds to multiplication under the Fourier transform, we have that $\|\varphi_j-\varphi_j^2\|_{A(H_j)}\leq \kappa_j^2$.  Recalling that $\lim_\cU \kappa_j^2=0$, it is in this sense that we see that each $\varphi_j$ is an almost idempotent element of the Banach algebra $A(H_j)$.

We need the following standard perturbation lemma, which is an immediate
consequence of \cite[Proposition~2.6]{DawsHorvath21}.
We also remark that in the discussion preceding that proposition, Daws and
Horv{\'a}th also refer to Johnson's
\cite[Lemma~2.1]{Johnson1988}
for an earlier statement of the idempotent perturbation result.
\begin{lem}\label{tar:banach-idempotent}
Let $\mathcal E$ be a unital complex Banach algebra with $\|1\|=1$.
Suppose $a\in\mathcal E$ satisfies
\[
 \|a\|\leq1,
 \qquad \varepsilon:=\|a-a^2\|\leq\frac1{16}.
\]
Then there is an idempotent $b\in\mathcal E$ such that $\|b-a\|\leq6\varepsilon.$
\end{lem}

\begin{proof}
By \cite[Proposition~2.6]{DawsHorvath21}, since
$\varepsilon<1/4$, there is an idempotent
$b\in\mathcal E$ such that
\[
\begin{aligned}
\|b-a\|
&\le \left(\|a\|+\frac12\right)
     \left((1-4\varepsilon)^{-1/2}-1\right)\\
&\le \frac32\cdot4\varepsilon
=6\varepsilon.
\end{aligned}
\]
Here we used $\|a\|\le1$, $4\varepsilon\le1/4$,
and the elementary inequality
$(1-t)^{-1/2}-1\le t$ for $0\le t\le1/4$.
\end{proof}

Now consider those indices $j$ for which $\kappa_j^2\leq \frac{1}{16}$.  For these indices, there are idempotents $b_j\in A(H_j)$ with $\|\varphi_j-b_j\|_{A(H_j)}\leq 6\kappa_j^2$.  It remains to show that $b_j=1_{\{e\}}$ for $\cU$-many such $j$.  We first note that idempotent elements of $A(H_j)$ take values in $\{0,1\}$.  

Since the $A(H_j)$-norm dominates the $\ell^{\infty}$ norm, 
we see that $|b_j(e)-\varphi_j(e)|\leq 6\kappa_j^2$.  Recalling that $\varphi_j(e)=\tr(\rho_j(e))=1$, we see that $b_j(e)=1$ provided that $6\kappa_j^2<1/2$.  It remains to show, for $\cU$-many $j$, that $b_j(g)=0$ for all
$g\in H_j\setminus\{e\}$. Set
\[
\delta_j:=
\max_{g\in H_j\setminus\{e\}}
|\operatorname{tr}_{D_j}(\rho_j(g))|,
\]
where the maximum is taken to be zero if $H_j=\{e\}$.
We claim that $\lim_{\mathcal U}\delta_j=0$.
If this were not the case, then there exists $\varepsilon>0$ and
$g_j\in H_j\setminus\{e\}$ for $\mathcal U$-many
$j$ such that
$|\operatorname{tr}_{D_j}(\rho_j(g_j))|
\ge\varepsilon$.
By \eqref{eq:rho-separation} (where we may take $s=1$),
$\|\rho_j(g_j)-1\|_{2,D_j}\ge1/2$
for $\mathcal U$-many $j$.
Thus $u:=(\rho_j(g_j))_{\mathcal U}$
is a nonidentity element of
$K_{\mathbf D,\Phi_2(V)}$ with
$|\tau(u)|\ge\varepsilon$, contradicting
the trace-zero property of that group.

Consequently, for $\cU$-many $j$ and all $g\in H_j\setminus \{e\}$, we have
\[
|b_j(g)|
\le |\varphi_j(g)|+6\kappa_j^2
\le \delta_j+6\kappa_j^2<1.
\]
Since $b_j(g)\in\{0,1\}$ and $b_j(e)=1$,
we conclude that $b_j=1_{\{e\}}$ for
$\mathcal U$-many $j$. Hence
\[
\|\varphi_j-1_{\{e\}}\|_{A(H_j)}
\le6\kappa_j^2\longrightarrow_{\mathcal U}0.
\]

This completes the proof of Proposition \ref{prop:nearidentity} and thus the proof of Theorem \ref{thm:special}.

\appendix
\section{The quantifier complexity of the sentence $\sigma$}\label{app:complexity}

In Subsection \ref{subsec:proofofmaintheorem}, we asserted that the maximum of the formulae (1)--(12) is
(equivalent to) a $\forall_5$-formula, so that existentially quantifying
over the parameters yields the $\exists_{6}$-sentence $\sigma$.  In this appendix, we show that the proof of Theorem~2.6 goes through with a
modified list of conditions, each of which is $\forall\exists$. Consequently, the
sentence $\sigma$ can be taken to be $\exists\forall\exists$. The key
observation is that the relative bicommutant $\A_{\mathbf 2}$ of $G_{\mathbf 4}$
is defined by a \emph{quantifier-free} formula: it is the relative commutant of
the two parameters $\Theta(b_1),\Theta(b_2)$.

Throughout this appendix:
\begin{itemize}
\item $\eta$ is the constant from Fact \ref{aff:external-expander}, and $L=1+8\sqrt2\,\eta$ is the
constant from Proposition \ref{aff:ambient-rounding}.
\item We set
\[
c:=2\eta,\qquad L_1:=1+4\sqrt2\,c,\qquad L_2:=\sqrt{1+c^2}.
\]
\end{itemize}

Let $M$ be a tracial von Neumann algebra and let $b_1,b_2,F,V\in U(M)$. Write
$\vec b=(b_1,b_2)$. We consider the following formulae, where $x$ ranges over
the unit ball $M_1$ and $u$ ranges over $U(M)$:
\begin{align*}
\alpha(x,\vec b)&:=\max_{i=1,2}\|[x,b_i]\|_2,\\
\kappa(u,V,F)&:=\|VuV^*-uFuF^*\|_2,\\
\psi(u)&:=\max\big(\alpha(u,\vec b),\kappa(u,V,F)\big),\\
\psi_c(x)&:=\sup_{u\in U(M)}\big(\|[x,u]\|_2\dotminus 2L_1\psi(u)\big),\\
\psi'(x)&:=\max\big(\psi_c(x),\alpha(x,\vec b)\big).
\end{align*}
Parameters are suppressed when no confusion can arise. Note that $\alpha$,
$\kappa$, and $\psi$ are quantifier-free. We also set
\[
A^{\vec b}:=\{b_1,b_2\}'\cap M,\qquad
K:=\{u\in U(A^{\vec b}):\kappa(u)=0\},\qquad
C:=K'\cap A^{\vec b}.
\]
Since $b_1,b_2$ are unitaries, $A^{\vec b}$ is a von Neumann subalgebra of $M$.
Moreover, we can identify the zerosets of $\alpha$ and $\psi$ as
\[
Z(\alpha)=(A^{\vec b})_1\quad\text{and}\quad Z(\psi)=K.
\]

As in the main text, we freely quantify over $U(M)$. This does not affect
quantifier complexity, because $U(M)$ is defined in $M_1$ by a quantifier-free
formula with a linear modulus.

\begin{lem}\label{lem:app-qf-bicommutant}
Let $\mathbf q=(q_j)$ be a sequence of prime powers with $q_j\ge3$ and
$\lim_{\mathcal U}q_j=\infty$, and set $n_j:=q_j(q_j-1)$. Let
$b_1,b_2\in U(\M_{\mathbf n})$ be as in Section \ref{sec:affine}, and set
$\A:=(G_{\mathbf n}'\cap \M_{\mathbf n})'\cap \M_{\mathbf n}$. Then
\begin{enumerate}
    \item $\A=\{b_1,b_2\}'\cap \M_{\mathbf n}$, and
    \item for all $x\in(\M_{\mathbf n})_1$, we have
$d(x,\A_1)\le 2\eta\,\alpha(x,\vec b)$.
\end{enumerate}

Consequently, in $\M_{\mathbf 4}$ we have
$\A_{\mathbf 2}=\{\Theta(b_1),\Theta(b_2)\}'\cap \M_{\mathbf 4}$ and
\[
d(x,(\A_{\mathbf 2})_1)\le c\,\alpha(x,\Theta(\vec b))
\qquad\text{for all }x\in(\M_{\mathbf 4})_1 .
\]
\end{lem}

\begin{proof}
Since
$\lambda_q(H_q)'\cap M_n(\mathbb C)=B_q$ and $B_q'\cap M_n(\mathbb C)=A_q$,
Lemma \ref{rec:internal-group} gives $\A=\prod_{\mathcal U}A_{q_j}$.

Fix $q$ and a contraction $x\in M_n(\mathbb C)$,
where $n=q(q-1)$. Let $p=p_q$ and set $y:=E_{A_q}(pxp)\in(A_q)_1$. Exactly as in
the proof of Lemma \ref{aff:selector-rounding}, we have:
\begin{itemize}
\item $p$ is central in $A_q$, $E_{A_q}$ restricts on $pM_n(\mathbb C)p$ to the
expectation onto $M_d(\mathbb C)\otimes1$, and
$[pxp,pb_{i,q}]=p[x,b_{i,q}]p$. Together with (4.7), this gives
\[
\|pxp-y\|_2\le 2\eta\max_{i=1,2}\|[x,b_{i,q}]\|_2 .
\]
\item $x-pxp$ is orthogonal to $pM_n(\mathbb C)p$ and
\[
\|x-pxp\|_2^2=\|(1-p)x\|_2^2+\|px(1-p)\|_2^2\le 2\tr_n(1-p)=2/q .
\]
\end{itemize}
Combining the two, we have
\[
d(x,(A_q)_1)\le\|x-y\|_2\le 2\eta\max_{i=1,2}\|[x,b_{i,q}]\|_2+\sqrt{2/q}.
\]

Given $x\in(\M_{\mathbf n})_1$, choose
contractive representatives $x=(x_j)_{\mathcal U}$, and let $y_j\in(A_{q_j})_1$
be as above. Then $y:=(y_j)_{\mathcal U}\in \A_1$, and since $\lim_\cU q_j=\infty$, we have $
\|x-y\|_2\le 2\eta\,\alpha(x,\vec b)$.

Since $b_i\in\prod_{\mathcal U}B_{q_j}$ and
$[A_{q_j},B_{q_j}]=0$, we have $\A\subseteq\{b_1,b_2\}'\cap \M_{\mathbf n}$.  The other containment is clear.

The final statement follows by applying
the isomorphism $\Theta$, which satisfies $\Theta(\A)=\A_{\mathbf 2}$.
\end{proof}

\begin{lem}\label{lem:app-definability}
Let $M$ be a tracial von Neumann algebra and let $b_1,b_2,F,V\in U(M)$. Write
$A:=A^{\vec b}$, and let $K$ and $C$ be as above. Suppose that:
\begin{enumerate}
\item[(a)] $d(x,A_1)\le c\,\alpha(x)$ for all $x\in M_1$, and
\item[(b)] $d(a,K)\le\kappa(a)$ for all $a\in U(A)$.
\end{enumerate}
Then:
\begin{enumerate}
\item $d(u,U(A))\le\sqrt2\,c\,\alpha(u)$ for all $u\in U(M)$.
\item $\psi$ is $L_1$-almost near with zeroset $K$.
\item If, moreover, $K$ is a subgroup of $U(M)$, then for all $x\in M_1$, we have:
\begin{itemize}
\item $\psi_c(x)=\sup_{k\in K}\|[x,k]\|_2$;
\item $Z(\psi')=C_1$;
\item $d(x,C_1)\le L_2\,\psi'(x)$.
\end{itemize}
\end{enumerate}
\end{lem}

\begin{proof}
(1) is proven just as in the proof of Lemma \ref{lem:specialdefinable} and only uses assumption (a).

For (2), it is clear that $Z(\psi)=K$. As in the proof of Lemma \ref{lem:specialdefinable}, $\kappa$ is 3-Lipschitz on
$U(M)$; this uses only that $V$ and $F$ are unitary. Let $u\in U(M)$ and take
$a\in U(A)$ as in the polar decomposition argument of part (1), whence $\|u-a\|_2\leq \sqrt{2}c\alpha(u)$. By (b),
\begin{align*}
d(u,K)&\le\|u-a\|_2+\kappa(a)\\
&\le 4\|u-a\|_2+\kappa(u)\\
&\le 4\sqrt2\,c\,\alpha(u)+\kappa(u)\\
&\le L_1\,\psi(u).
\end{align*}

(3) Since $1\in K$ and $y\mapsto\|[x,y]\|_2$ is 2-Lipschitz and bounded by 2 for
$x\in M_1$, Fact \ref{fact:quant} and (2) give $\psi_c(x)=\sup_{k\in K}\|[x,k]\|_2$. Hence
$\psi'(x)=0$ if and only if $x\in A$ and $x$ commutes with $K$, proving that
$Z(\psi')=C_1$.

Now fix $x\in M_1$ and set $a:=E_A(x)\in A_1$. Then
$\|x-a\|_2\le d(x,A_1)\le c\,\alpha(x)$. Applying Lemma \ref{syn:averaging} inside $A$ to the
subgroup $K\le U(A)$, and using that $E_A$ is $A$-bimodular and contractive, yields
\[
\|a-E_C(a)\|_2
\le\sup_{k\in K}\|[a,k]\|_2
=\sup_{k\in K}\|E_A([x,k])\|_2
\le\psi_c(x).
\]
Since $x-a\perp a-E_C(a)$ and $E_C(a)=E_C(x)\in C_1$, we get
\[
d(x,C_1)^2\le\|x-a\|_2^2+\|a-E_C(a)\|_2^2\le c^2\alpha(x)^2+\psi_c(x)^2
\le L_2^2\,\psi'(x)^2. \qedhere
\]
\end{proof}

Now consider parameter variables $b_1,b_2,W,F,V$ ranging over unitaries. Let (C1)-(C7)
be the formulae (1)--(7) of Subsection \ref{subsec:proofofmaintheorem}, with the parameters
$\Theta(\vec b)$, $\Theta(W)$, $\Flip_{\mathbf 2}$ replaced by the
variables $\vec b$, $W$, $F$. Let $G$ denote the zeroset of
$\varphi_{\mathrm{aff}}(\cdot,\vec b,W,F)$, and let
$\varphi_{\mathrm{aff},c}$ be as in Proposition \ref{prop:definability}. We replace the formulae (8)-(12) of Subsection \ref{subsec:proofofmaintheorem} with the following
formulae, written first informally and then officially via Fact \ref{fact:quant}:
\begin{enumerate}
\item[(C8)] $R_c(\alpha)$, that is,
\[
\sup_{x\in M_1}\inf_{y\in M_1}
\max\big(\alpha(y),\|x-y\|_2\dotminus c\,\alpha(x)\big).
\]
\item[(C9)] $\sup_{x\in Z(\varphi_{\mathrm{aff},c}),\ y\in Z(\alpha)}\|[x,y]\|_2$, that is,
\[
\sup_{x,y\in M_1}\Big(\big(\|[x,y]\|_2\dotminus 2c\,\alpha(y)\big)
\dotminus 2\varphi_{\mathrm{aff},c}(x)\Big).
\]
\item[(C10)] $R_{L_1}(\psi)$.
\item[(C11)] $\sup_{u\in U(A^{\vec b})}\big(d(u,K)\dotminus\kappa(u)\big)$, that is,
\[
\sup_{u}\Big(\inf_{v}\min\big(2,\|u-v\|_2+L_1\psi(v)\big)
\dotminus\kappa(u)\dotminus 4\sqrt2\,c\,\alpha(u)\Big).
\]
\item[(C12)] $\sup_{u\in K}\|u^2-1\|_2$, that is,
$\sup_{u}\big(\|u^2-1\|_2\dotminus 2L_1\psi(u)\big)$.
\item[(C13)] $\sup_{x,y\in C_1}\|[x,y]\|_2$, that is,
\[
\sup_{x,y\in M_1}\big(\|[x,y]\|_2\dotminus 2L_2\psi'(x)\dotminus 2L_2\psi'(y)\big).
\]
\end{enumerate}
In (C11) and (C12), $u,v$ range over unitaries. The constants come from:
\begin{itemize}
\item (C9): Proposition \ref{prop:definability}(1), and condition (a) of
Lemma~\ref{lem:app-definability} with Lipschitz constant 2;
\item (C11): Lemma~\ref{lem:app-definability}(1),(2), using that
$u\mapsto d(u,K)\dotminus\kappa(u)$ is 4-Lipschitz;
\item (C12): Lemma~\ref{lem:app-definability}(2), using that $u\mapsto u^2$ is
2-Lipschitz;
\item (C13): Lemma~\ref{lem:app-definability}(3).
\end{itemize}

It is now straightforward to verify that each of (C1)-(C13) is equivalent to a $\forall\exists$-formula. Consequently,
\[
\sigma:=\inf_{b_1,b_2,W,F,V}\ \max\big((\mathrm{C}1),\dots,(\mathrm{C}13)\big)
\]
is an $\exists\forall\exists$-sentence.

\begin{lem}\label{lem:app-holds}
$\sigma^{\M_{\mathbf 4}}=0$, as witnessed by the parameters
$\Theta(b_1)$, $\Theta(b_2)$, $\Theta(W)$, $\Flip_{\mathbf 2}$, $V$, where $V$ is a special
unitary of $\M_{\mathbf 4}$.
\end{lem}

\begin{proof}
(C1)--(C7) vanish for the reasons given in Subsection \ref{subsec:proofofmaintheorem}. By
Lemma~\ref{lem:app-qf-bicommutant}, we have $A^{\Theta(\vec b)}=\A_{\mathbf 2}$
and $\alpha$ is $c$-almost near. Hence (C8) vanishes. Since
$(\A_{\mathbf 2})'\cap \M_{\mathbf 4}=\B_{\mathbf 2}=G_{\mathbf 4}'\cap \M_{\mathbf 4}$,
(C9) vanishes.

For \(u\in \mathrm{U}(\A_{\mathbf2})\) we have
$u\Flip_{\mathbf 2}u\Flip_{\mathbf 2}=u\otimes u$. Therefore
\[
\kappa(u)=\|VuV^*-u\otimes u\|_2,\qquad K=K_{\mathbf 2,V},\qquad C=C_{\mathbf 2,V}.
\]
The three defining properties of a special unitary now give the remaining
conditions:
\begin{itemize}
\item Property (1) gives assumption (b)
of Lemma~\ref{lem:app-definability}, so part (2) of that lemma shows that
(C10) vanishes.  Fact \ref{fact:quant} then identifies (C11) with property (1).
\item Property (2) gives that (C12) vanishes.
\item $K_{\mathbf 2,V}$ is a group, so Lemma~\ref{lem:app-definability}(3)
applies; together with property (3), (C13) vanishes.\qedhere
\end{itemize}
\end{proof}

\begin{prop}\label{prop:app-thm26}
Suppose that $\M_{\mathbf 4}\equiv \M_{\mathbf r}$. Then there is a sequence
$\mathbf t$ with $\mathbf r/\mathbf t^2=1$ such that $\M_{\mathbf t^2}$ has a
special unitary.
\end{prop}

\begin{proof}
By Lemma~\ref{lem:app-holds}, $\sigma^{\M_{\mathbf r}}=0$. Since $\M_{\mathbf r}$ is
countably saturated, we may take witnesses
$\tilde b_1,\tilde b_2,\widetilde W,\widetilde F,\widetilde V\in U(\M_{\mathbf r})$
to the vanishing of (C1)-(C13). For ease of notation, we drop the tildes from the notation for
$\alpha,\kappa,\psi,\psi'$, and write $G:=G_{\mathbf r}$,
$A^\#:=A^{\tilde{\vec b}}$, and $K,C$ for the corresponding sets.

As in Subsection \ref{subsec:proofofmaintheorem}, (C1)-(C7) show that
$\varphi_{\mathrm{aff}}(\cdot,\tilde{\vec b},\widetilde W,\widetilde F)$ is
$L$-almost near, $\widetilde F^2=1$, and $G$ satisfies the hypotheses of
Theorem \ref{thm:frame}. By Proposition \ref{prop:definability}, the zeroset of $\varphi_{\mathrm{aff},c}$
is $(G'\cap \M_{\mathbf r})_1$, and
$d(x,(G'\cap \M_{\mathbf r})_1)\le\varphi_{\mathrm{aff},c}(x)$ for all
$x\in(\M_{\mathbf r})_1$.

By (C8) and countable
saturation, $\alpha$ is $c$-almost near, that is, $A^\#$ satisfies condition (a)
of Lemma~\ref{lem:app-definability}. In particular, by Fact \ref{fact:quant}, the official
formula in (C9) computes
$\sup\{\|[x,y]\|_2:x\in(G'\cap \M_{\mathbf r})_1,\,y\in A^\#_1\}$.  

We claim that $A^\#=(G'\cap M_{\mathbf r})'\cap \M_{\mathbf r}$.  Indeed, $\varphi_{\mathrm{aff}}(u)\ge\|[u,\tilde b_i]\|_2$, so
$G\subseteq A^\#$. Since each \(\tilde b_i\) commutes with \(G\), we have
\(\tilde b_i\in G'\cap\M_{\mathbf r}\) for \(i=1,2\). Hence
\[
 (G'\cap\M_{\mathbf r})'\cap\M_{\mathbf r}
 \subseteq \{\tilde b_1,\tilde b_2\}'\cap\M_{\mathbf r}=A^\#.
\]
On the other hand, (C9) shows that $A^\#$ commutes with $G'\cap \M_{\mathbf r}$, so
$A^\#\subseteq(G'\cap \M_{\mathbf r})'\cap \M_{\mathbf r}$.

We next show that $K$ is a group.  By (C6), $\widetilde FG\widetilde F\subseteq G'$.
Using $\widetilde F=\widetilde F^*$ and taking relative commutants in
$\M_{\mathbf r}$, we get
$A^\#=G''\subseteq(\widetilde FG\widetilde F)'=\widetilde F G'\widetilde F$, so
\[
\widetilde FA^\#\widetilde F\subseteq G'=(A^\#)'.
\]
Thus $[A^\#,\widetilde FA^\#\widetilde F]=0$. It follows that
$u\mapsto u\widetilde Fu\widetilde F$ is a group homomorphism on $U(A^\#)$
commuting with adjoints. As $\Ad\widetilde V$ is an automorphism, $K$ is a
subgroup of $U(A^\#)$, as desired.

We next check that the special unitary conditions hold for $K$ and $V$.  By (C10), $\psi$ is $L_1$-almost near. By Lemma~\ref{lem:app-definability}(1)
and Fact \ref{fact:quant}, (C11) says that condition (b) of
Lemma~\ref{lem:app-definability} holds.
(C12) then says that $K$ has exponent $2$.
By the previous paragraph and Lemma~\ref{lem:app-definability}(3), $C$ is defined by
$\psi'$ with constant $L_2$, so (C13) says that $C$ is abelian.

Let
$\Pi:\M_{\mathbf r}\to \M_{\mathbf t^2}$ be the isomorphism given by
Theorem \ref{thm:frame}, with $\mathbf r/\mathbf t^2=1$. By the above and the conclusions of
Theorem \ref{thm:frame}, $\Pi(A^\#)=\A_{\mathbf t}$ and
$\Pi\circ\Ad\widetilde F=\Ad(\Flip_{\mathbf t})\circ\Pi$. Hence, for
$u\in U(A^\#)$, we have
\[
\Pi\big(u\widetilde Fu\widetilde F^*\big)=\Pi(u)\otimes\Pi(u),
\]
so that $\kappa(u)=\|\Pi(\widetilde V)\Pi(u)\Pi(\widetilde V)^*-\Pi(u)\otimes\Pi(u)\|_2$.
Consequently, $\Pi(K)=K_{\mathbf t,\Pi(\widetilde V)}$ and 
$\Pi(C)=C_{\mathbf t,\Pi(\widetilde V)}$. 
By the previous paragraph, we have that $\Pi(\widetilde V)$ is a special unitary of
$\M_{\mathbf t^2}$, as desired.
\end{proof}

The rest of the argument
in Subsection \ref{subsec:proofofmaintheorem} is unchanged. In particular, $\sigma$ is an
$\exists\forall\exists$-sentence such that $\sigma^{\M_{\mathbf 4}}=0$, while
$\sigma^{(\M_{\mathbf 4})_s}>0$ for $s\notin 4^{\mathbb Z}$.
\bibliographystyle{siam}
\bibliography{reference}

\end{document}